\documentclass[11pt]{article}
\usepackage{amsmath,latexsym,epsf,mathabx}
\usepackage{makecell}
\usepackage{changepage}
\usepackage[all]{xy}
\usepackage{hyperref}
\usepackage{color}
\usepackage{ulem}

\usepackage{epsfig}
\usepackage{tikz}
\usepackage{graphicx}
\usepackage{caption}
\usepackage{subcaption}
\usepackage{indentfirst}
\usepackage{extarrows}
\usepackage{dsfont}

\begin{document}

\newcommand{\nc}{\newcommand}
\def\PP#1#2#3{{\mathrm{Pres}}^{#1}_{#2}{#3}\setcounter{equation}{0}}
\def\ns{$n$-star}\setcounter{equation}{0}
\def\nt{$n$-tilting}\setcounter{equation}{0}
\def\Ht#1#2#3{{{\mathrm{Hom}}_{#1}({#2},{#3})}\setcounter{equation}{0}}
\def\qp#1{{${(#1)}$-quasi-projective}\setcounter{equation}{0}}
\def\mr#1{{{\mathrm{#1}}}\setcounter{equation}{0}}
\def\mc#1{{{\mathcal{#1}}}\setcounter{equation}{0}}
\def\HD{\mr{Hom}}
\def\HC{\mr{Hom}_{\mc{C}}}
\def\AdT{\mr{Add}_{\mc{T}}}
\def\adT{\mr{add}_{\mc{T}}}
\def\Kb{\mc{K}^b(\mr{Proj}R)}
\def\kb{\mc{K}^b(\mc{P}_R)}
\def\AdpC{\mr{Adp}_{\mc{C}}}
\def\AdpD{\mr{Adp}_{\mc{D}}}
\newtheorem{Th}{Theorem}[section]
\newtheorem{Def}[Th]{Definition}
\newtheorem{Lem}[Th]{Lemma}
\newtheorem{Pro}[Th]{Proposition}
\newtheorem{Cor}[Th]{Corollary}
\newtheorem{Rem}[Th]{Remark}
\newtheorem{Exm}[Th]{Example}
\newtheorem{Sc}[Th]{}
\def\Pf#1{{\noindent\bf Proof}.\setcounter{equation}{0}}
\def\>#1{{ $\Rightarrow$ }\setcounter{equation}{0}}
\def\<>#1{{ $\Leftrightarrow$ }\setcounter{equation}{0}}
\def\bskip#1{{ \vskip 20pt }\setcounter{equation}{0}}
\def\sskip#1{{ \vskip 5pt }\setcounter{equation}{0}}
\def\mskip#1{{ \vskip 10pt }\setcounter{equation}{0}}
\def\bg#1{\begin{#1}\setcounter{equation}{0}}
\def\ed#1{\end{#1}\setcounter{equation}{0}}
\def\KET{T^{^F\bot}\setcounter{equation}{0}}
\def\KEC{C^{\bot}\setcounter{equation}{0}}

\def\jze{{ \begin{pmatrix} 0 & 0 \\ 1 & 0 \end{pmatrix}}\setcounter{equation}{0}}
\def\hjz#1#2{{ \begin{pmatrix} {#1} & {#2} \end{pmatrix}}\setcounter{equation}{0}}
\def\ljz#1#2{{  \begin{pmatrix} {#1} \\ {#2} \end{pmatrix}}\setcounter{equation}{0}}
\def\jz#1#2#3#4{{  \begin{pmatrix} {#1} & {#2} \\ {#3} & {#4} \end{pmatrix}}\setcounter{equation}{0}}

\title{ \bf Strongly $n$-AIR-tilting modules
\footnotetext{Version of \today}
\footnotetext{E-mail:~weiqingcao@jsnu.edu.cn (Weiqing Cao),~weijiaqun@njnu.edu.cn (Jiaqun Wei)}
\footnotetext{Weiqing cao was supported by the Science Foundation of Jiangsu Normal University (No. 21XFRS024). Jiaqun wei was supported by  the National Natural Science Foundation of China (Grant No. 11771212)}}

\smallskip
\author{\small Weiqing Cao$^{a}$, Jiaqun Wei$^{ b}$ \\
 \small    $^{a}$ School of Mathematics and Statistics, Jiangsu Normal University,
    Xuzhou 221116,  P.R.China\\
\small $^{b}$ School of Mathematics Science, Nanjing Normal University, Nanjing 210024,
P.R. China}
\date{}
\maketitle
\baselineskip 15pt
%
%
\begin{abstract}
\vskip 10pt%
%
%
We introduce the notion of (strongly) $n$-AIR-tilting modules, which is a high dimesion version of support $\tau$-tilting modules. The relations between them and $n$-silting modules and $n$-quasi-tilting modules, as well as generalized two-term silting complexes, are investigated. Our results particularly suggest a way to negate  the rank question for silting complexes.

\mskip\




\end{abstract}
%
\vskip 20pt

\section{Introduction}
Tilting theory is an important tool in the representation theory of artin algebras. Classical tilting modules were introduced by Brenner and Butler\cite{BB}. Happel\cite{Happel} and Miyashita \cite{Miyashita} initiated the study of tilting modules of finite projective dimension. Further generalizations of classical tilting modules, such as $\ast$-modules, infinitely generated tilting modules, quasi-tilting modules, etc. were presented in \cite{C,CT,AC,CDT, Wei4} etc. 

Adachi, Iyama and Reiten\cite{AIR} introduced the concept of support $\tau$-tilting modules as a generalization of classical tilting modules, and they showed that the mutations of support $\tau$-tilting modules are always possible. In particular, it was shown that  support $\tau$-tilting modules are in bijection with two-term silting complexes. 
Angeleri-H\"{u}gel, Marks and Vit\'{o}ria\cite{AMV} introduced silting modules over an arbitrary algebra as an infinitely generated version of support $\tau$-tilting modules. 
Mao\cite{MAO} introduced the concept of $n$-silting modules as a high dimension version of silting modules.

In this paper, we introduce a new high dimension version of support $\tau$-tilting modules (and also silting modules), which we call (strongly) $n$-AIR-tilting modules, and study their properties and characterizations. We compare them with other related concepts and investigate their relations with general silting complexes.

Our results can be described as the following diagram.

$$\xymatrix{&&   *+[F]{\thead{\text{strongly} \\ \text{$n$-quasi-tilting} \\ \text{modules}}}\ar[dr]^{\subset}\ar@{-}[d]^{\neq~~~~\neq}\ar@{-}@/^2.5pc/[dd]\\
*+[F]{{\thead{\text{$n$-tilting} \\ \text{modules}}} }\ar[r]^{\subset~~~} &  *+[F]{\thead{\text{strongly} \\ \text{ $n$-AIR-tilting} \\ \text{modules}}}\ar[ur]^{~~~\subset}\ar[r]^{~~~\subseteq}\ar[dr]^{\subseteq}&  *+[F]{\thead{\text{$n$-silting} \\ \text{modules}}}\ar[r]^{\subset}&  *+[F]{\thead{\text{$n$-quasi-tilting} \\ \text{modules}}}\\
&  *+[F]{\thead{\text{generalized two-term} \\ \text{silting complexes} \\ \text{ of length $n+1$}}}\ar[u]^{\Phi}&   *+[F]{\thead{\text{ $n$-AIR-tilting} \\ \text{modules}}}\ar[ur]^{\subset}\ar[l]_{~~\Psi}\\ 
} $$

\vskip 20pt
We note that the map $\Psi$ is established under the assumption that the rank question of silting complexes has a positive answer. In the case, the maps $\Phi, \Psi$ are inverse bijections, in particular, strongly $n$-AIR-tilting modules and $n$-AIR-tilting modules coincide with each other. Thus, if the two modules are different in general, then the answer to the
rank question for silting complexes is negative. By this meaning, a way to negate the rank question is to show the difference between strongly $n$-AIR-tilting modules and $n$-AIR-tilting modules.

The paper is organized as follows.
In Section 2 and Section 3, we introduce some special classes of modules and discuss their elementary properties.
In Section 4, we introduce the new concept of $n$-AIR-tilting modules and provide some characterizations. In particular, we describe the relations between (strongly) $n$-AIR-tilting modules, (strongly) $n$-quasi-tilting modules and $n$-silting modules. The inclusions in the above diagram are established in this section.
In Section 5, we introduce the generalized  two-term silting complexes of length $(n+1)$ and study their relations with (strongly) $n$-AIR-tilting modules. In particular, the maps $\Phi, \Psi$ are provided in this section.  Some interesting examples and questions are listed in Section 6.


\vskip 20pt

\section{Preliminaries}

Let $A$ be  an Artin algebra, mod$A$  the category of finitely generated right $A$-modules, $\mc{P}_{A}$ the category of finitely generated projective right $A$-modules.

Let $\mc{C}$ be a subcategory of $\mr{mod}A$ and $k\ge 1$. We denote by $\mr{Pres}^k(\mc{C})$ the {\it $k$-images} of $\mc{C}$, i.e., the modules $M$ such that there is an exact sequence $C_k\to \cdots\to C_1\to M\to 0$, where the exact sequence will be called a {\it $k$-$\mc{C}$-presentation} of $M$ and the correspondent truncated complex $C_k\to \cdots\to C_1\to 0$ will be called the a {\it truncated  $k$-$\mc{C}$-presentation} of $M$. A $k$-$\mc{C}$-presentation $C_k\to \cdots\to C_1\to M\to 0~(\spadesuit)$ is called an {\it approximate $k$-$\mc{C}$-presentation} if $\mathrm{Hom}(C,\spadesuit)$ is eaxct for each $C\in\mathcal{C}$, i.e., the induced sequence $\mathrm{Hom}(C,C_k)\to \cdots\to \mathrm{Hom}(C,C_1)\to \mathrm{Hom}(C,M)\to 0$ is exact. We denote $\mr{Appres}^{k}(\mathcal{C})$ the class of all modules which has an approximate $k$-$\mc{C}$-presentation.

%

Dually, we denote by $\mr{Copres}^k(\mc{C})$ the $k$-subs of $\mc{C}$, i.e., the modules $M$ such that there is an exact sequence $0\to M\to  C_1\to \cdots\to C_k$ which we call a {\it $k$-$\mc{C}$-copresentation}. 

In case $k=1$, $\mr{Pres}^1(\mc{C}) = \mr{Gen}(\mc{C}) = \mr{Fac}(\mc{C})$, where last two are usual notations which are denoted to be the class of quotients of modules in $\mc{C}$, and $\mr{Copres}^1(\mc{C}) = \mr{Cogen}(\mc{C}) = \mr{Sub}(\mc{C})$, where last two are usual notations which are denoted to be the class of submodules of modules in $\mc{C}$. It is clear that $\mr{Pres}^k(\mc{C})\subseteq \mr{Pres}^m(\mc{C})$ and that $\mr{Copres}^k(\mc{C})\subseteq \mr{Copres}^m(\mc{C})$ for any $k>m\ge 1$.

For a subcategory $\mc{C}$ and a subclass $\mc{M}$ of morphisms, we say that, respectively,
\begin{verse}
\noindent (1) $\mc{C}$ is closed under extensions if for any exact sequence $0\to L\to M\to N\to 0$ with $L,N\in\mc{C}$ it holds that $M\in\mc{C}$.

\noindent (2) $\mc{C}$ is closed under cokernels of $\mc{M}$-morphisms if for any exact sequence $L\stackrel{f}{\to} M\to N\to 0$ with $f\in\mc{M}$ it holds that $N\in\mc{C}$.

\noindent (3) $\mc{C}$ is closed under kernels of $\mc{M}$-morphisms if for any exact sequence $0\to L\to M\stackrel{f}{\to} N\to 0$ with $f\in\mc{M}$ it holds that $L\in\mc{C}$.

\noindent (4) $\mc{C}$ is closed under $k$-images if $\mr{Pres}^k(\mc{C})\subseteq \mc{C}$.

\noindent (5) $\mc{C}$ is closed under $k$-subs if $\mr{Copres}^k(\mc{C})\subseteq \mc{C}$.
\end{verse}

\vskip 10pt 
 Let $T\in \mr{mod}A$.  We fix the following  two exact sequences

          \centerline{$\sigma_{_{T}}: P_{n}\to P_{n-1}\to\cdots\to P_{0}\to T\to 0$,}
            \noindent and

            \centerline{$\sigma_{_{A}}: A\to T_{0}\to T_{1}\to \cdots\to T_{n}\to 0$}
\noindent with $P_{i}\in\mc{P}_{A}$ and $T_{i}\in\mr{add}(T)$  for $0\leq i \leq n$. Such $\sigma_{_T}$ will be called an $(n+1)$-projective presentation of $T$. We denote $\sigma^{\bullet}_{_{T}}$ the brutal truncated complex of $\sigma_{_{T}}$ at the term $T$, i.e.,

\centerline{$\sigma^{\bullet}_{_{T}}: 0\to P_{n}\to P_{n-1}\to\cdots\to P_{0}\to 0$.}

\vskip 10pt 
Let $1 \leq k \leq n$. We consider the following classes of $A$-modules
\begin{verse}
 $\mc{D}_{\sigma_{_{T}}}\ =\{M \in \mr{mod}A ~| ~\mr{~the~indueced~sequence~} \mr{Hom}_{A}(P_{0},M)\to\mr{Hom}_{A}(P_{1},M)\to\cdots\to\mr{Hom}_{A}(P_{n},M)\to 0$
  is exact $\}$.

 $\mc{D}_{\sigma_{_{T}},k}\ =\{M \in \mr{mod}A ~| ~\mr{~the~indueced~sequence~} \mr{Hom}_{A}(P_{n-k},M) \to \mr{Hom}_{A}(P_{n-k+1},M) \to \cdots \to \mr{Hom}_{A}(P_{n},M) \to 0$
  is exact $\}$.

 $\mc{D}_{\sigma_{_{A}}}\ =\{M \in \mr{mod}A ~| ~\mr{~the~indueced~sequence~} \mr{Hom}_{A}(T_{n},M)\to\mr{Hom}_{A}(T_{n-1},M)\to\cdots\to\mr{Hom}_{A}(T_0,M)\to\mr{Hom}_{A}(A,M)\to 0$
  is exact $\}$.

 $\mc{D}_{\sigma_{_{A}},k}\ =\{M \in \mr{mod}A ~| ~\mr{~the~indueced~sequence~} \mr{Hom}_{A}(T_{k-1},M)\to \mr{Hom}_{A}(T_{k-2},M)\to\cdots\to\mr{Hom}_{A}(T_0,M)\to \mr{Hom}_{A}(A,M)\to0$
  is exact $\}$.

\end{verse}

In particular, $\mc{D}_{\sigma_{_T},n} = \mc{D}_{\sigma_{_T}}$.
It is easy to see that $\mc{D}_{\sigma_{_T},k+1}\subseteq \mc{D}_{\sigma_{_T},k}$.

The unbounded derived (respectively, homotopy) category of mod$A$ will be denoted by $D$($A$) (resp., $K(A)$), or just $D$ (resp., $K$). If we restrict ourselves to bounded complexes, we use the usual superscript $b$.
For a subcategory $\mc{X}$ of the derived category of $A$-module category $D(A)$, we denote by $\mc{X}^{\perp_{>0}}$ the  subcategory consisting of the  object $Y$ in $D(A)$ such that $\mr{Hom}_{D}(X,Y[i])=0$ for all $i>0$ and all $X\in\mc{X}$. If the subcategory consists of a single object $X$, we write just $X^{\perp_{>0}}$.
We denote by ${D}^{\leq0}\subseteq({D}{(A}))$ the complexes whose homologies are concentrated on non-positive terms.

Let $\mathrm{F}_i: \mathcal{C}\to\mathcal{D}$ be a set of functors with $i\in I$. We denote $$\mathrm{KerF}_{i\in I}: = \bigcap_{i\in I} \mathrm{KerF}_i = \{X\in \mathcal{C}\ |\ \mathrm{F}_i(X)=0 \mathrm{\ for\ all\ } i\in I\}.$$

For instance, for an $A$-module $T$, $$\mathrm{KerExt}_A^{1\le i\le n}(T,-) := \{X\in \mathrm{mod}A\ |\ \mathrm{Ext}_A^i(T,X)=0 \mathrm{\ for\ all\ } 1\le i\le n\}.$$

In case that the subcategory $\mathcal{C}=\mr{add}(T)$ for some object $T$, we simply replace $\mathcal{C}$ with $T$ in the correspondent notations.

\vskip 10pt

\section{Properties of $\mc{D}_{\sigma_{_T}}$, $\mc{D}_{\sigma_{_T},k}$, $\mc{D}_{\sigma_{_A}}$}

Now, we provide properties of the classes $\mc{D}_{\sigma_{_T}}$, $\mc{D}_{\sigma_{_T},k}$,  $\mc{D}_{\sigma_{_A}}$.

\bg{Lem}\label{1}
Let $T\in\mr{mod}A$ and $\sigma_{_{T}}$ be an   $(n+1)$-projective presentation of $T$. The class $\mc{D}_{\sigma_{_{T}}}\subseteq \mr{KerExt}^{i}_{A}(T,-)$ for $1\leq i \leq n$. And the class $\mc{D}_{\sigma_{_{T}},k}\subseteq\mr{KerExt}_A^{i}(T,-)$ for all $n-k+1\leq i \leq n$, where $1\leq k \leq n$.
\ed{Lem}
\Pf. It is easy to prove by the definitions of the functors $\mr{KerExt}^{i}_{A}(T,-)$.
\ \hfill $\Box$

\bg{Lem}\label{n-image}
Let $T\in\mr{mod}A$ and $\sigma_{_{T}}$ be an $(n+1)$-projective presentation of $T$. Then

$(1)$ For any $n\ge k\ge 1$, $\mc{D}_{\sigma_{_{_{T}}},k}$ is closed under $k$-images, i.e.,  $\mr{Pres}^{k}(\mc{D}_{\sigma_{_{T}},k})\subseteq\mc{D}_{\sigma_{_{T}},k}$.  In particular, $\mc{D}_{\sigma_{_{T}}}$ is closed under $n$-images.

$(2)$ If $T\in\mc{D}_{\sigma_{_{T}}}$, then $\mr{Pres}^{k}(T)\subseteq\mc{D}_{\sigma_{_{T}},k}\subseteq\mr{KerExt}_A^{n-k+1\leq i \leq n}(T,-)$, where $1\leq k \leq n$. In particular, $\mr{Pres}^{n}(T)\subseteq \mc{D}_{\sigma_{_{T}}}$ in the case.

$(3)$  $\mc{D}_{\sigma_{_{T}}}$ is closed under extensions.

$(4)$ $\mc{D}_{\sigma_{_{T}}}$ is closed under cokernels of monomorphisms.

$(5)$ $\mc{D}_{\sigma_{_{T}}}$ is closed under kernels of $\mr{Hom}_{A}(T,-)$-epic  epimorphisms.
\ed{Lem}
\Pf.
(1) We proceed induction on $k$.

Take any $M\in\mr{Pres}^{1}(\mc{D}_{\sigma_{_{T}},1})$. Then there is an epimorphism $C_{1}\stackrel{g_{1}}\to M\to 0$, where  $C_{1}\in\mc{D}_{\sigma_{_{T}},1}$.  Take any $f:P_{n}\to M$ and consider the following diagram.

 $$
 \xymatrix{
 & &P_{n}\ar[r]^{f_{n}}\ar[d]_{f}\ar@{.>}[dl]_{g}                &  P_{n-1}\ar[r] \ar@{.>}[dll]^{h}                 &  P_{n-2}\ar[r]&\cdots                      \\
 &C_{1}\ar[r]_{g_{1}}&M  \ar[r]                      &0 &  & \\
}
$$

As $P_n$ is projective, there is a homomorphism $g:P_{n}\to C_{1}$ such that $f=g_{1}g$. Since $C_{1}\in\mc{D}_{\sigma_{_{T}},1}$, there exists $h:P_{n-1}\to C_{1}$ such that $g=hf_{n}$. Hence we have $f=g_{1}hf_{n}$. So, $M\in\mc{D}_{\sigma_{_{T}},1}$.



Now assume $k>1$ and the result holds for all $i<k$. Take $M\in\mr{Pres}^{k}(\mc{D}_{\sigma_{_{T}},k})$. Then there is an exact sequence $0\to M_1\stackrel{g_{2}}\to C_{1}\stackrel{g_{1}}\to M\to 0$, where $C_1\in\mc{D}_{\sigma_{_{T}},k}$ and $M_1\in \mr{Pres}^{k-1}(\mc{D}_{\sigma_{_{T}},k})$.   We only need to prove that $\mr{Hom}_{A}(P_{n-k},M)\to \mr{Hom}_{A}(P_{n-k+1},M)\to \mr{Hom}_{A}(P_{n-k+2},M) $ is exact by the induction assumptions. This is shown by considering the following diagram.

 $$
 \xymatrix{
 &&&P_{n-k+3}\ar[r]^{f_{n-k+3}}&P_{n-k+2}\ar[r]^{f_{n-k+2}} \ar@{.>}[dl]_{h_{2}}                &  P_{n-k+1} \ar[d]|{h}\ar[r]^{f_{n-k+1}}\ar@{.>}[dl]|{h_{1}}  \ar@{.>}[dll]^{s_{2}}                &  P_{n-k}\ar[r]\ar@{.>}[dll]_{s_{1}}&\cdots                       \\
 &&0\ar[r]&M_1\ar[r]^{g_{2}}&C_{1}\ar[r]^{g_{1} }                       & M\ar[r] &  0&  \\
}
$$

Take  $h\in\mr{Hom}_{A}(P_{n-k+1},M)$ satisfying $hf_{n-k+2}=0$.
Since $P_{n-k+1}$ is projective, there is a map $h_{1}:P_{n-k+1}\to C_{1}$ such that $h=g_{1}h_{1}$. It is easy to see that there is some $h_2:P_{n-k+2}\to M_1$ such that $h_1f_{n-k+2}=g_{2}h_{2}$. It follows $g_2h_2f_{n-k+3}=0$ and then $h_2f_{n-k+3}=0$ since $g_2$ is monic. As $M_1\in \mr{Pres}^{k-1}(\mc{D}_{\sigma_{_{T}},k}) \subseteq\mr{Pres}^{k-1}(\mc{D}_{\sigma_{_{T}},k-1})$, there exists $s_{2}:P_{n-k+1}\to M_1$ such that $s_{2}f_{n-k+2}=h_{2}$. Thus, $(h_{1}-g_{2}s_{2})f_{n-k+2}=0$ and then there exists $s_{1}:P_{n-k}\to C_{1}$ such that $s_{1}f_{n-k+1}=h_{1}-g_{2}s_{2}$ since $C_1\in\mc{D}_{\sigma_{_{T}},k}$. It follows that $h=g_{1}h_{1}=g_{1}(s_{1}f_{n-k+1}+g_{2}s_{2})=g_{1}s_{1}f_{n-k+1}$. So, $M\in\mc{D}_{\sigma_{_{T}},k}$.

When $k=n$, we obtain that  $\mc{D}_{\sigma_{_{T}}}=\mc{D}_{\sigma_{_{T}},n}$ is closed under $n$-image.

(2) It follows from (1) and Lemma \ref{1} (1).

(3)
Let $0\to L\stackrel{\alpha}\to M\stackrel{\beta}\to N\to 0$ be exact with $L,N\in\mc{D}_{\sigma_{_{T}}}$. Take any $f\in\mr{Hom}_{A}(P_{n},M)$. Since $N\in\mc{D}_{\sigma_{_{T}}}$, there exists $g: P_{n-1}\to N$ such that $gf_{n}=\beta f$. Since $P_{n-1}$ is projective, there exists $h:P_{n-1}\to M$ such that $g=\beta h$. Then $\beta(f-hf_{n})=\beta f-g f_{n}=0$; hence there is some $l: P_{n}\to L$ such that $f-hf_{n}=\alpha l$, and hence $f=hf_{n}+\alpha l$. Since $L\in\mc{D}_{\sigma_{_{T}}}$, there exists $s:P_{n-1}\to L$ such that $l=s f_{n}$. Then $f=hf_{n}+\alpha sf_{n}=(h+\alpha s)f_{n} $; hence $\mr{Hom}_{A}(P_{n-1},M)\to\mr{Hom}_{A}(P_{n},M)\to 0$ is exact.
 $$
 \xymatrix{
 & &P_{n}\ar[r]^{f_{n}}\ar[d]|{f}\ar@{.>}[dl]_{l}                &  P_{n-1}\ar[r]\ar@{.>}[d]^{g}   \ar@{.>}[dl]|{h}    \ar@{.>}[dll]|{s}             &  P_{n-2}\ar[r]&\cdots                      \\
 0\ar[r]&L\ar[r]_{\alpha}&M  \ar[r]_{\beta}                      &N \ar[r]&0  & \\
}
$$

Let $\gamma\in\mr{Hom}_{A}(P_{i},M) $ such that $\gamma f_{i+1}=0$, where $1 \leq i \leq n-1$. Then there exists $\delta: P_{i-1}\to N$ such that $\delta f_{i}=\beta \gamma$. Since $P_{i-1}$ is projective,  there is a map $\epsilon: P_{i-1}\to M$ such that $\delta=\beta \epsilon$. Then $\beta(\gamma -\epsilon f_{i})=\beta \gamma-\delta f_{i}=0$; hence there is some $\zeta: P_{i}\to L$ such that $\alpha\zeta=\gamma-\varepsilon f_{i}$. Then  $\alpha\zeta f_{i+1}=0$. Since $\alpha$ is injective, $\zeta  f_{i+1}=0$. Then there is a map  $\eta:P_{i-1}\to L$ such that $\zeta=\eta f_{i}$ because $L\in\mc{D}_{\sigma_{_{T}}}$.  Hence $\gamma=\alpha\eta f_{i}+\varepsilon f_{i}=(\alpha\eta+\varepsilon)f_{i}$. Then
$\mr{Hom}_{A}(P_{i-1},M)\to\mr{Hom}_{A}(P_{i},M)\to\mr{Hom}_{A}(P_{i+1},M)$ is exact.

 $$
 \xymatrix{
 &\cdots\ar[r] &P_{i+1}\ar[r]^{f_{i+1}}                &  P_{i} \ar[d]_{\gamma}\ar[r]^{f_{i}}\ar@{.>}[dl]|{\zeta}                 &  P_{i-1}\ar[r]\ar@{.>}[dll]^{\eta}\ar@{.>}[d]^{\delta}\ar@{.>}[dl]^{\epsilon}&\cdots                       \\
 &0\ar[r]&L\ar[r]_{\alpha }                       & M\ar[r]_{\beta} &  N\ar[r]&0  \\
}
$$

(4)
Let $0\to L\stackrel{\alpha}\to M\stackrel{\beta}\to N\to 0$ be exact with $L,M\in\mc{D}_{\sigma_{_{T}}}$.
Take $f\in\mr{Hom}_{A}(P_{n},N)$. Since $P_{n}$ is projective, there exists $g:P_{n}\to M$ such that $f=\beta g$. Since $M\in\mc{D}_{\sigma_{_{T}}}$, there is a map $h:P_{n-1}\to M$ such that $g=hf_{n}$. Then $f=\beta h f_{n}$.
 $$
 \xymatrix{
 && &P_{n}\ar[r]^{f_n}\ar[d]_{f}\ar@{.>}[dl]_{g}                &  P_{n-1}\ar[r] \ar@{.>}[dll]^{h}                 &  P_{n-2}\ar[r]&\cdots                      \\
0\ar[r]& L\ar[r]_{\alpha}&M\ar[r]_{\beta}&N  \ar[r]                      &0 &  & \\
}
$$
Take $\gamma\in\mr{Hom}_{A}(P_{n-1}, N)$ satisfying $\gamma f_{n}=0$. Since $P_{n-1}$ is projective, there is a map $\delta:P_{n-1}\to M$ such that $\gamma=\beta \delta$. Since $\beta \delta f_{n}=0$, there exists $\epsilon:P_{n}\to L$ such that $\alpha \epsilon=\delta f_{n}$.  Then there is a map $\eta:P_{n-1}\to L$ such that $\epsilon= \eta f_{n}$ because $L\in \mc{D}_{\sigma_{_{T}}}$. Since $(\delta-\alpha \eta)f_{n}=0$ and $M\in\mc{D}_{\sigma_{_{T}}}$, there exists $\theta:P_{n-2}\to M$ such that $\delta=\alpha\eta+\theta f_{n-1}$. Then $\gamma=\beta \delta= \beta\theta f_{n-1}$.
 $$
 \xymatrix{
 & &P_{n}\ar[r]^{f_{n}}\ar@{.>}[dl]|{\epsilon}               &  P_{n-1}\ar[r]^{f_{n-1}}\ar[d]^{\gamma}   \ar@{.>}[dl]|{\delta}    \ar@{.>}[dll]|{\eta}             &  P_{n-2}\ar[r]\ar@{.>}[dll]^{\theta}&\cdots                      \\
 0\ar[r]&L\ar[r]_{\alpha}&M  \ar[r]_{\beta}                      &N \ar[r]&0  & \\
}
$$
Take $a\in\mr{Hom}_{A}(P_{i}, N)$ satisfying $af_{i+1}=0$. Since $P_{i}$ is projective, there is a map $b:P_{i}\to M$ such that $a=\beta b$. Since $\beta b f_{i+1}=0$, there exists $c:P_{i+1}\to L$ such that $\alpha c=b f_{i+1}$.
 Since $\alpha c f_{i+2}=0$ and  $\alpha$ is injective, we obtain $cf_{i+2}=0$. Since $L\in\mc{D}_{\sigma_{_{T}}}$,
there is a map $d:P_i \to L$ such that $c= d f_{i+1}$ because $L\in \mc{D}_{\sigma_{_T}}$. Since $(b-\alpha d)f_{i+1}=0$ and $M\in\mc{D}_{\sigma_{_{T}}}$, there exists $d:P_{i-1}\to M$ such that $b=\alpha d+a f_{i}$. Then $a=\beta b= \beta a f_{i}$.
$$
 \xymatrix{
 \cdots\ar[r]& P_{i+2}\ar[r]^{f_{i+2}} &P_{i+1}\ar[r]^{f_{i+1}}\ar@{.>}[dl]|{c}               &  P_{i}\ar[r]^{f_{i}}\ar[d]^{a}   \ar@{.>}[dl]|{b}    \ar@{.>}[dll]|{d}             &  P_{i-1}\ar[r]\ar@{.>}[dll]^{d}&\cdots                      \\
 0\ar[r]&L\ar[r]_{\alpha}&M  \ar[r]_{\beta}                      &N \ar[r]&0  & \\
}
$$

(5)
Let $0\to L\stackrel{\alpha}\to M\stackrel{\beta}\to N\to 0$ be a $\mr{Hom}_{A}(T,-)$-exact exact sequence   with $M,N\in\mc{D}_{\sigma_{_{T}}}$. Take any $f\in\mr{Hom}_{A}(P_{n},M)$.  Since $M\in\mc{D}_{\sigma_{_{T}}}$, there is a map $s:P_{n-1}\to M$ such that $sf_{n}=\alpha f$. Since $\beta s f_{n}=\beta\alpha f=0$, there exists $s^{'}:P_{n-2}\to N$ such that $s^{'}f_{n-1}=\beta s$. Since $P_{n-2}$ is projective, there is a map $t:P_{n-2}\to M$ such that $s^{'}=\beta t$. Then $\beta(s-tf_{n-1})=0$; hence there is some $t^{'}:P_{n-1}\to L$ such that $s=\alpha t^{'}+tf_{n-1}$. Since $\alpha(f-t^{'}f_{n})=0 $ and $\alpha$ is injective, we obtain $f=t^{'}f_{n}$.

$$
 \xymatrix{
 & &P_{n}\ar[r]^{f_{n}} \ar[d]^{f}               &  P_{n-1} \ar@{.>}[d]_{s}\ar[r]^{f_{n-1}}\ar@{.>}[dl]|{t^{'}}                 &  P_{n-2}\ar[r]\ar@{.>}[dl]^{t}\ar@{.>}[d]^{s^{'}}&\cdots                       \\
 &0\ar[r]&L\ar[r]_{\alpha }                       & M\ar[r]_{\beta} &  N\ar[r]&0  \\
}
$$

Take $a\in\mr{Hom}_{A}(P_{i},L)$ satisfying $af_{i+1}=0$, where $2\leq i \leq n-1$. Then there are $b:P_{i-1}\to M$ and $c:P_{i-2}\to N$ such that $bf_{i}=\alpha a $ and $cf_{i-1}=\beta b$.  Since $P_{i-2}$ is projective, there exists $d:P_{i-2}\to M$ such that $c=\beta d$.  Since $\beta(b-d f_{i-1})=0$, there is a map $e:P_{i-1}\to L$ such that $b=d f_{i-1}+\alpha e$. Then $\alpha a= bf_i=\alpha e f_{i}$. Since $\alpha $ is injective, $a=e f_{i}$. So, $\mr{Hom}_{A}(P_{i-1},L)\to\mr{Hom}_{A}(P_{i},L)\to \mr{Hom}_{A}(P_{i+1},L) $ is exact.
$$
 \xymatrix{
\cdots\ar[r] & P_{i+1}\ar[r]^{f_{i+1}}&P_{i}\ar[r]^{f_{i}} \ar[d]^{a}               &  P_{i-1} \ar@{.>}[d]_{b}\ar[r]^{f_{i-1}}\ar@{.>}[dl]|{e}                 &  P_{i-2}\ar[r]\ar@{.>}[dl]^{d}\ar@{.>}[d]^{c}&\cdots                       \\
 &0\ar[r]&L\ar[r]_{\alpha }                       & M\ar[r]_{\beta} &  N\ar[r]&0  \\
}
$$

It remains to prove that  $\mr{Hom}_{A}(P_{0},L)\to \mr{Hom}_{A}(P_{1},L)\to\mr{Hom}_{A}(P_{2},L)$ is exact. This proof is similar to the above because  $0\to L\to M\to N\to 0$ is Hom$_{A}(T,-)$-exact.

\ \hfill $\Box$

\vskip 10pt
Similarly, we have the following result.

\bg{Cor}\label{k}
Let $T\in\mr{mod}A$ and $\sigma_{_{T}}$ be an $(n+1)$-projective presentation of $T$. Then

$(1)$  $\mc{D}_{\sigma_{_{T}},k}$ is closed under extension.

$(2)$ $\mc{D}_{\sigma_{_{T}},k}$ is closed under cokernels of monomorphisms.

\ed{Cor}

\vskip 10pt
Now we turn to study the properties of $\mc{D}_{\sigma_{_{A}}}$. Recall that $\sigma_{_{A}}$ is the exact sequence $A\to T_{0}\to T_{1}\to \cdots\to T_{n}\to 0$, see Section 2.

\bg{Lem}\label{DA}
$(1)$ $\mc{D}_{\sigma_{_{A}}}$ is closed under kernels of  $\mr{Hom}_{A}(T,-)$-epic morphisms.

$(2)$ ${D}_{\sigma_{_{A}},k}\subseteq\mr{Appres}^{1}(T)$ for any $1\le k\le n$. If moreover $T\in\mc{D}_{\sigma_{_{A}}}$, then $\mc{D}_{\sigma_{_{A}},k}=\mr{Appres}^{k}(T)$.
\ed{Lem} 

\Pf.
(1) Let $0\to L\stackrel{\alpha}\to\ M\stackrel{\beta} \to N$ be an  exact sequence with $M,N\in{\mc{D}}_{\sigma_{_{A}}}$ such that Hom$_{A}(T,\beta)$ is surjective.  Take any $f\in\mr{Hom}_{A}(A,L)$. Since $M\in{\mc{D}}_{\sigma_{_{A}}}$, there is $g:T_{0}\to M$ such that $gf_{0}=\alpha f$. Then there is $h:T_{1}\to N$ such that $hf_{1}=\beta g$. Since Hom$_{A}(T_{1},\beta)$ is surjective, there is $s:T_{1}\to M$ such that $h=\beta s$. Since $\beta(g-sf_{1})=0$, we have $t:T_{0}\to L$ with $g=\alpha t+ sg$. Then $\alpha f=gf_0=(\alpha t+ sg)f_{0}=\alpha tf_{0}$. Since $\alpha$ is injective, $f=t f_{0}$.
$$
 \xymatrix{
 & &A\ar[r]^{f_{0}} \ar[d]^{f}               &  T_{0} \ar@{.>}[d]_{g}\ar[r]^{f_{1}}\ar@{.>}[dl]|{t}                 &  T_{1}\ar[r]\ar@{.>}[dl]^{s}\ar@{.>}[d]^{h}&\cdots\ar[r]  & T_{n}\ar[r]&0                     \\
 & 0\ar[r]                     & L\ar[r]_{\alpha} &M\ar[r]_{\beta} &N &&\\
}
$$

Now take any $i$ such that $0\leq i \leq n-2$. Take $l\in\mr{Hom}_{A}(T_{i},L)$ satisfying $lf_{i}=0$. Then there are $m:T_{i+1}\to M$ and $n:T_{i+2}\to N$ such that $\alpha l=mf_{i+1}$ and $nf_{i+2}=\beta m$. Since Hom$_{A}(T,\beta)$ is surjective, there is $s_{2}:T_{i+2}\to M$ such that $n=\beta s_{2}$.  Since $\beta(m-s_{2}f_{i+2})=0$, there is some $s_{1}:T_{i+1}\to M$ such that $m=s_{2}f_{i+2}+\alpha s_{1}$. Then $\alpha l=mf_{i+1}=\alpha s_{1} f_{i+1}$. Since $\alpha$ is injective, $l=s_{1}f_{i+1}$. Hence $\mr{Hom}_{A}(T_{i+1},L)\to\mr{Hom}_{A}(T_{i},L)\to \mr{Hom}_{A}(T_{i-1},L)\to 0$ is exact.

 $$
 \xymatrix{
\cdots\ar[r] &T_{i-1}\ar[r]^{f_{i}}              &  T_{i} \ar[d]_{l}\ar[r]^{f_{i+1}}                &  T_{i+1}\ar[r]^{f_{i+2}}\ar@{.>}[dl]^{s_{1}}\ar@{.>}[d]^{m}&T_{i+2}\ar[r]\ar@{.>}[d]^{n}\ar@{.>}[dl]^{s_{2}}&\cdots\ar[r]  & T_{n}\ar[r]&0                     \\
                     & 0\ar[r]&L\ar[r]_{\alpha} &  M\ar[r]_{\beta}  & N&&\\
}
$$

It follows that $L\in {\mc{D}}_{\sigma_{_{A}}}$. Thus  ${\mc{D}}_{\sigma_{_{A}}}$ is closed under kernels of Hom$_{A}(T,-)$-epic morphisms.

(2)  Take any $M\in{\mc{D}}_{\sigma_{_{A}},k}$. Then there is an integer $m$ such that there exists an epimorphism $f
:A^{m}\to M$. Now consider the following diagram where the first row is $\sigma_{_{A}}^{m}$.

$$
 \xymatrix{
&A^{m}\ar[r]\ar[d]^{f}& T^{m}_{0}\ar[r]\ar@{.>}[dl]^{g}     &  T_{1}^{m} \ar[r]            &  \cdots\ar[r] &T_{n}^{m} \ar[r]&0&&                      \\
 &M                     &  && & & && \\
}
$$

Since $M\in{\mc{D}}_{\sigma_{_{A}},k}$, there is $g: T_0^m\to M$ which makes the above diagram commutative. It follows that $g$ is also surjective. This shows that ${\mc{D}}_{\sigma_{_{A}},k}\subseteq \mr{Gen}(T)$.

Now we prove that  ${D}_{\sigma_{_{A}},k}\subseteq\mr{Appres}^{k}(T)$ in case that $T\in\mc{D}_{\sigma_{_{A}}}$.

Take any $M\in{\mc{D}}_{\sigma_{_{A}},k}$, then there is a  $\mr{Hom}_{A}(T,-)-$epic morphism, $0\to M_{1}\to T_{0}^{'}\to M\to 0$, where $T_{0}^{'}\in \mr{add}(T)$. By (1) and what we just proved, we have that $M_{1}\in{\mc{D}}_{\sigma_{_{A}},k}\subseteq\mr{Gen}(T)$. Now by repeating this process to $M_{1}$ and continuously,   we  obtain that  $M\in \mr{Appres}^{\infty}(T)$. In particular, $\mc{D}_{\sigma_{_{A}},k}\subseteq\mr{Appres}^{k}(T)$ follows.

Conversely, we proceed it by induction on $k$. 

Let $k=1$. Take $M\in\mr{Appres}^{1}(T)$ and any $f\in\mr{Hom}_{A}(A,M)$. Then there is a Hom$_{A}(T,-)$-epic epimorphism  $T_1^{'}\stackrel{\alpha_{1}}\to M$, where $T_1^{'}\in\mr{add}(T)$.  

 $$
 \xymatrix{
 & &A\ar[r]^{g_{0}}\ar[d]_{f}\ar@{.>}[dl]_{g}                &  T_{0}\ar[r] \ar@{.>}[dll]^{h}                 &  T_{1}\ar[r]&\cdots                      \\
 &T_1^{'}\ar[r]_{\alpha_{1}}&M  \ar[r]                      &0 &  & \\
}
$$

Since $M\in{\mc{D}}_{\sigma_{_{A}},k}$, there is $g: T_0^m\to M$ which makes the above diagram commutative. It follows that $g$ is also surjective. This shows that ${\mc{D}}_{\sigma_{_{A}},k}\subseteq \mr{Gen}(T)$. Hence there is a  $\mr{Hom}_{A}(T,-)-$epic morphism, $0\to M_{1}\to T_{0}^{'}\to M\to 0$, where $T_{0}^{'}\in \mr{add}(T)$. Then, $M_{1}\in{\mc{D}}_{\sigma_{_{A}},k}\subseteq\mr{Gen}(T)$. Repeating this process to $M_{1}$ and continuously,   we  could obtain  ${D}_{\sigma_{_{A}},k}\subseteq\mr{Appres}^{k}(T)$.
Hence there is a map $g:A\to T_1^{'}$ such that $f=\alpha_{1}g$. Since $T_1^{'}\in{\mc{D}}_{\sigma_{_{A}},1}$, there exists $h:T_{0}\to T_1^{'}$ such that $g=hg_{0}$. Hence we have $f=\alpha_{1}hg_{0}$. So, $M\in{\mc{D}}_{\sigma_{_{A}},1}$ and $\mr{Appres}^{1}(T)\subseteq  {\mc{D}}_{\sigma_{_{A}},1}$.

Let $k>1$. Take $M\in\mr{Appres}^{k}(T)$. Then there is a Hom$_{A}(T,-)$-exact exact sequence $ T_k^{'}\stackrel{\alpha_{k}}\to \cdots\stackrel{\alpha_{2}}\to T_1^{'}\stackrel{\alpha_{1}}\to M\to 0$, where $T_i^{'}\in\mr{add}(T)$ for $1 \leq i \leq k$. We only need to prove that $\mr{Hom}_{A}(T_{k-1},M)\to \mr{Hom}_{A}(T_{k-2},M)\to \mr{Hom}_{A}(T_{k-3},M) $ is exact. Consider the following diagram.

 $$
 \xymatrix{
 &A\ar[r]^{g_{0}}\ar@{.>}[dl]_{h_{k}}& T_{0}\ar[r]^{g_{1}}\ar@{.>}[dll]^{s_{k}}\ar@{.>}[dl]^{h_{k-1}}&\cdots\ar[r]&T_{k-3}\ar[r]^{g_{k-2}} \ar@{.>}[dl]_{h_{2}}                &T_{k-2}\ar[d]|{h}\ar[r]^{g_{k}}\ar@{.>}[dl]|{h_{1}}  \ar@{.>}[dll]^{s_{2}}                &  T_{k-1}\ar[r]\ar@{.>}[dll]_{s_{1}}&\cdots                       \\
 T_{k}^{'}\ar[r]_{\alpha_{k}}\ar[r]&T_{k-1}^{'}\ar[r]&\cdots\ar[r]&T_{2}^{'}\ar[r]_{\alpha_{2}}&T_{1}^{'}\ar[r]_{\alpha_{1} }                       & M\ar[r] &  0&  \\
}
$$
Take any $h\in\mr{Hom}_{A}(T_{k-2},M)$ satisfying $hg_{k-2}=0$.
Since $\alpha_{1}$ is Hom$_{A}(T,-)$-epic, there is a map $h_{1}:T_{k-2}\to T_1^{'}$ such that $h=\alpha_{1}h_{1}$. It is easy to see that there are $h_{i}:T_{k-1-i}\to T_i^{'}$ such that $h_{i-1}g_{k-i}=\alpha_{i}h_{i}$ for $2\leq i\leq k$, where $A:=T_{-1}$. Since $T_k^{'}\in {\mc{D}}_{\sigma_{_{A}},k}$, there exists $s_{k}:T_{0}\to T_k^{'}$ such that $s_{k}g_{0}=h_{k}$. Since $(h_{k-1}-\alpha_{k}s_{k})g_{0}=0$ and $T_{k-1}^{'}\in {\mc{D}}_{\sigma_{_{A}},k}$, there exists $s_{k-1}:T_{1}\to T_{k-1}^{'}$ such that $\alpha_{k}s_{k}+s_{k-1}g_{1}=h_{k-1}$. Similarly, since $T_{i}^{'}\in{\mc{D}}_{\sigma_{_{A}},k}$, there are $s_{i}:T_{k-i}\to T_{i}^{'}$ such that $h_{i}=s_{i}g_{k-i}+\alpha_{i+1}s_{i+1}$ for $1\leq i \leq k-1$. Then $h=\alpha_{1}h_{1}=\alpha_{1}(s_{1}g_{k}+\alpha_{2}s_{2})=\alpha_{1}s_{1}g_{k}$. So, $M\in {\mc{D}}_{\sigma_{_{A}},k}$  and $\mr{Appres}^{k}(T)\subseteq  {\mc{D}}_{\sigma_{_{A}},k}$.

As a conclusion, we have $\mr{Appres}^{k}(T)=  {\mc{D}}_{\sigma_{_{A}},k}$.

\ \hfill $\Box$

\section{$n$-AIR-tilting modules}
We introduce $n$-AIR-tilting modules and study the related propertied in this section. Then, we prove our main results Theorem \ref{AIR-snqt} and Theorem \ref{main2}. We give  relationships of related  notions in this section.

We introduce the following definitions, where AIR is the first letters of the last names of the authors of \cite{AIR}.

\bg{Def}\rm{\label{AIR}}
Let $T\in\mr{mod}A$. 

 $(1)$ We call $T$ an $n\text{-}AIR\text{-}tilting$ module if it satisfies the following conditions:

     $(i)$ there exists an exact sequence $\sigma_{_{T}}: P_{n}\to P_{n-1}\to \cdots\to P_{0}\to T\to 0$ with $P_{i}\in\mc{P}_{A}$ for $0 \leq i \leq n$ such that $0\to \mr{Hom}_{A}(T,T)\to\mr{Hom}_{A}(P_{0},T)\to \cdots \to\mr{Hom}_{A}(P_{n},T)\to 0$ is exact, that is $T\in\mc{D}_{\sigma_{_{T}}}$;

     $(ii)$ there exists an exact sequence $\sigma_{_{A}}:A\to T_{0}\to T_{1}\to\cdots\to T_{n}\to 0$  with $T_{i}\in\mr{add}(T)$ for $0\leq i \leq n$ such that $0\to \mr{Hom}_{A}(T_{n},T)\to\mr{Hom}_{A}(T_{n-1},T)\to \cdots \to\mr{Hom}_{A}(A,T)\to 0$ is exact, that is $T\in\mc{D}_{\sigma_{_{A}}}$.

 In the case, we also say that $T$ is $n$-AIR-tilting with respect to $(\sigma_{_{T}},\sigma_{_{A}})$.

$(2)$ We call $T$ a $strongly~n\text{-}AIR\text{-}tilting$ module if it is $n$-AIR-tilting with respect to $(\sigma_{_{T}},\sigma_{_{A}})$ and satisfies $\mc{D}_{\sigma_{_{T}}}=\mc{D}_{\sigma_{_{A}}}$.
\ed{Def}

When $n=1$, we could see  that $1$-AIR-tilting modules coincide with support $\tau$-tilting modules \cite{AIR}.

Recall that an $A$-module $T$ is $n\text{-}tilting$ \cite{Happel}\cite{Miyashita} if it satisfies the following conditions:
\begin{verse}
    
$(1)$ The projective dimension of  $T$ is at most $n$;

$(2)$ $\text{Ext}^{i}_{A}(T,T)=0$ for all $i>0$;

$(3)$ There is a finite coresolution of $A$ by objects of $\text{add}(T)$,  i.e., an exact sequence $0\to A\to T_{0}\to\cdots\to T_{n}\to 0$, with each $T_{i}\in\text{add}(T)$, for some $n$.
\end{verse}

Modules satisfying the first two of the above conditions are called {\it pretilting modules}.

Note that, if $T$ is an $n$-tilting module, then the module ${_B}T$ is also $n$-tilting, where $B=\mathrm{End}T_A$, and it holds that $T_A$ is faithful balanced, i.e., $A\simeq \mathrm{End}_BT$ \cite{Miyashita}. Sometime, we also say that $T$ is an $n$-tilting $B$-$A$-bimodule.

\vskip 10pt
It is easy to see that $n$-tilting modules are $n$-AIR-tilting modules. Our first result on $n$-AIR-tilting modules provides a close relation between $n$-AIR-tilting modules and $n$-tilting modules.

We use $\mr{ann}(T)$ to denote the annihilator ideal of $T$.

\bg{Pro}\label{AIR-tilting}
Let $T\in\mr{mod}A$ and be an $n$-AIR-tilting module. Then $T$ is an $n$-tilting $B$-$\bar{A}$-bimodule, where $\bar{A}=A/\mr{ann}(T)$ and $B=\mr{End}T_{A}$.
\ed{Pro}
\Pf.
Since $T$ is an  $n$-AIR-tilting module, there is a Hom$_{A}(-,T)$-exact exact sequence
$$P_{n}\to P_{n-1}\to \cdots\to P_{0}\to T\to 0,~~~~~~~~~~~~~~~~~~~~~~~~~~~~~~~~~~~~~~~~~~~~~~~~~~~(\ast_{1})$$
 where $P_{i}\in\mc{P}_A$ for $0 \leq i \leq n$.
Then we obtain the induced exact sequence
$$0\to\mr{Hom}_{A}(T,T)\to\mr{Hom}_{A}(P_{0},T)\to\cdots\to\mr{Hom}_{A}(P_{n},T)\to 0.~~~~~~~~~~~~~~~(\ast_{2})$$
Since $B=\mr{Hom}_{A}(T,T)$ and $\mr{Hom}_{A}(P_{i},T)\in\mr{add}(_{B}T)$,  the above exact sequence gives an exact sequence
$0\to B\to T^{'}_{0}\to T^{'}_{1}\to \cdots\to T^{'}_{n}\to 0$ where $T^{'}_{i}\in\mr{add}({_{B}T})$.

Moreover, there is a $\mr{Hom}_{A}(-,T)$-exact exact sequence $A\to T_{0}\to T_1\to\cdots\to T_{n}\to 0$ by Definition \ref{AIR}, where $T_{i}\in\mr{add}(T_{A})$ for $0\leq i\leq n$. Then we have the induced exact sequence
$$0\to \mr{Hom}_{A}(T_{n},T)\to\mr{Hom}_{A}(T_{n-1},T)\to\cdots\to \mr{Hom}_{A}(T_{0},T)\to \mr{Hom}_{A}(A,T)\to 0.$$

It gives an exact sequence $$0\to P^{'}_n\to P^{'}_{n-1}\to\cdots\to P^{'}_0\to T\to 0,~~~~~~~~~(\ast_{3})$$

where $P^{'}_{i}:=\mr{Hom}_{A}(T_{i},T)\in\mc{P}_B$ for $0\leq i\leq n$.

Applying Hom$_{B}(-,T)$ to $(\ast_{3})$, we obtain the following commutative diagram
$$
 \xymatrix{
&A\ar[r]\ar[d]^{\lambda}& T_{0}\ar[r]\ar[d]     &  T_{1} \ar[r]\ar[d]              &  \cdots\ar[r] &T_{n} \ar[r]\ar[d]&0&&                      \\
0\ar[r] &\mr{Hom}_{B}({T},{T})\ar[r]                     & \mr{Hom}_{B}(P^{'}_{0},T)\ar[r] &\mr{Hom}_{B}(P^{'}_{1},T)\ar[r]& \cdots\ar[r]& \mr{Hom}_{B}(P^{'}_{n},T)\ar[r]&0 && \\
}
$$
Since $T_{i}\simeq \mr{Hom}_{B}(P^{'}_{i},T)$ naturally, the morphism $\lambda:A\to \mr{Hom}_{B}(T,T)$ is surjective   and the second row is an exact sequence. Thus, Ext$^{i}_{B}(T,T)=0$ for $i>0$. It follows that $T$ is  an $n$-tilting $B$-module.

As $\lambda:A\to \mr{Hom}_{B}(T,T) = \mr{End}{_{B}T}$ is surjective, we have the exact sequence
$0\to \mr{ann} (T)\to A\stackrel{\lambda}\to \mr{End}{_{B}T}\to 0$, by the definition of the morphism $\lambda$. Then $\mr{End}{_{B}T}\simeq \bar{A}$ and hence $T$ is an $n$-tilting ${\bar{A}}$-module. It follows that $T$ is an $n$-tilting $B$-$\bar{A}$-bimodule.
\ \hfill $\Box$

\vskip 10pt
Recall that a module $T$  is faithful if $\mathrm{ann}(T)=0$.
As a corollary of the above result, we see that an $n$-tilting module is exactly a faithful $n$-AIR-tilting module. The following result provides another proof of this fact.

\bg{Pro}
Let $T\in\mr{mod}A$ and $T$ be an  $n$-AIR-tilting module. Assume $A\in\mr{Cogen}(T)$.
Then $T$ is an $n$-tilting module.
\ed{Pro}
\Pf. Since $T$ is an $n$-AIR-tilting module, there are two Hom$_{A}(-,T)$-exact exact sequences $\sigma_{_{T}}: P_{n}\stackrel{f_{n}}\to P_{n-1}\to\cdots\to P_{0}\to T\to 0$ and  $\sigma_{_A}: A\stackrel{g_{0}}\to T_{0}\to T_{1}\to \cdots\to T_{n}\to 0$, where $P_{i}\in\mc{P}_A$, $T_{j}\in\mr{add}(T)$  for $0\leq i,j \leq n$.
Since $A\in\mr{Cogen}T$, there are injective morphisms $f:P_{n}\to T^{'}$ and  $s:A\to T^{''}$, where  $T^{'},T^{''}\in\mr{add}(T)$. Then we have the following  commutative diagrams, since $T\in \mathcal{D}_{\sigma_{_T}} \bigcap\mathcal{D}_{\sigma_{_A}}$.

$$
 \xymatrix{
 & &P_{n}\ar[r]^{f_{n}} \ar[d]^{f}               &  P_{n-1} \ar[r]^{f_{n-1}}\ar@{.>}[dl]|{g}                 &  P_{n-2}\ar[r]&\cdots \ar[r]&T\ar[r]&0                      \\
 &&   T^{'}               && &  & &  \\
}
$$
$$
 \xymatrix{
 & &A\ar[r]^{g_{0}} \ar[d]^{s}               &  T_{0} \ar[r]\ar@{.>}[dl]|{t}                 & T_{1}\ar[r]&\cdots   \ar[r]&T_{n}\ar[r]&0                    \\
 &&   T^{''}                  &  & &  &&\\
}
$$

It follows that $f_{n}:P_{n}\to P_{n-1}$ and
$g_{0}:A\to T_{0}$ are injective.  Then  pd$T\leq n$ and there is an exact sequence $0
\to A\to T_{0}\to T_{1}\to \cdots\to T_{n}\to 0$.
Thus, it is easy to see that  $T$ is an $n$-tilting module.
\ \hfill $\Box$

\vskip 10pt
We recall the following definition about generalizations of quasi-tilting modules.

\bg{Def}\rm{\label{qtdef}}
Let $T\in\mr{mod}A$. 

$(1)$ We say that $T$ is $n\text{-}quasi\text{-}tilting$ if $\mr{Pres}^{n}(T)=\mr{Pres}^{n+1}(T)\subseteq\mr{KerExt}^{i}_{A}(T,-)$  for $1 \leq i\leq n.$ 

$(2)$ We say that $T$ is $strongly~{n}\text{-}quasi\text{-}tilting$ if $\mr{Pres}^{n}(T)=\mr{Pres}^{n+1}(T)$ and  $\mr{Pres}^{k}(T)\subseteq \mr{KerExt}^{n-k+1\leq i \leq n}_{A}(T,-)$, for each $1\leq k \leq n$. 
\ed{Def}

It is easy to see that if $T$ is an $n$-quasi-tilting module, then $\mr{Pres}^{n}(T)=\mr{Appres}^{\infty}(T)$.

\vskip 10pt
The following result gives a relation between $n$-AIR-tilting modules and $n$-quasi-tilting modules.

\bg{Pro}\label{AIR-nqt}
Let $T\in\mr{mod}A$ and $T$  be an  $n$-AIR-tilting module. Then $T$ is an $n$-quasi-tilting module.
\ed{Pro}
\Pf. As usual, let $\bar{A}=A/\mr{ann}(T)$.  It is clear that the inclusion $\mr{mod}\bar{A}\hookrightarrow \mr{mod}A$ is fully faithful, and $\mr{mod}\bar{A}$ can be identified with the subcategory $\{M\in\mr{mod}A\ |\ M\cdot \mr{ann}(T)=0\}$. In particular, the subcategory $\mr{Gen}(T)$, and hence $\mr{Pres}^n(T)$, is in $\mr{mod}\bar{A}$. By Proposition \ref{AIR-tilting}, $T$ is an $n$-tilting $\bar{A}$-module, it follows that Pres$^{n}_{\bar{A}}(T)=\mr{Pres}^{n+1}_{\bar{A}}(T)$ by \cite{Wei4}. Then the above arguments imply that Pres$^{n}_{A}(T)=\mr{Pres}^{n+1}_{A}(T)$. By Lemma \ref{n-image}, we obtain that Pres$^{n}_{A}(T)\subseteq \mc{D}_{\sigma_{_{T}}}\subseteq \mr{KerExt}^{i}_{A}(T,-)$ for $1\leq i\leq n$. It follows that $T$ is an $n$-quasi-tilting  module.
\ \hfill $\Box$

\vskip 10pt
We introduce the following definition.

\bg{Def}\rm{
Let $T\in\mr{mod}A$. We say that $T$ is of  $ann\text{-}faithful ~dimension ~at ~least ~n$ (briefly, $ann\text{-}faith.dimT\geq n$) if there is an exact sequence $A\to T_{1}\to\cdots \to T_{n}$  such that the induced sequence $\mr{Hom}_{A}(T_{n},T)\to\cdots\to\mr{Hom}_{A}(T_{1},T)\to \mr{Hom}_{A}(A,T)\to 0$ is exact, where $T_{i}\in\mr{add}(T)$ for $1 \leq i \leq n$. }
\ed{Def}

\noindent {\bf{Remark}} For any  $T\in\mr{mod}A$, $ann\text{-}faith.dimT\geq 1$ always holds since there is a morphism $f_{1}: A\to T_{1}$ such that $\mr{Hom}_{A}(f_{1},T)$ is surjective, for some $T_{1}\in\mr{add}(T)$.

\vskip 10pt
We now give a characterization of $n$-AIR-tilting modules in term of strongly $n$-quasi-tilting modules.

\bg{Th}\label{AIR-snqt}
Let $T\in\mr{mod}A$. Then the following are equivalent.

 $(1)$ $T$ is an $n$-AIR-tilting module.

 $(2)$  $T$ is a strongly $n$-quasi-tilting module such that ann-faith.dim$T\geq n$.
\ed{Th}

\Pf.(1)$\Rightarrow $(2) By the definition of $n$-AIR-tilting modules, we obtain ann-faith.dim$T\geq n$. By Lemma \ref{n-image}$(2)$ and Proposition \ref{AIR-nqt}, $T$ is a strongly $n$-quasi-tilting module.

(2)$\Rightarrow$(1) Take a minimal $(n+1)$-projective presentation $P_{n}\stackrel{f_{n}}\to P_{n-1}\to \cdots\to P_{0}\stackrel{f_{0}}\to T\to 0$.
Since Gen$(T)\subseteq\mr{KerExt}_{A}^{n}(T,-)\simeq \mr{KerExt}_{A}^{1}(\mr{Coker}(f_n),-)$, we obtain that Hom$_{A}(f_{n},T)$ is surjective. Since $T\in\mr{Pres}^{n}(T)\subseteq\mr{KerExt}^{1\leq i \leq n}_{A}(T,-)$, the induced sequence $\mr{Hom}_{A}(P_{0},T)\to\mr{Hom}_{A}(P_{1},T)\to\cdots\to\mr{Hom}_{A}(P_{n},T)$ is exact. So the condition $(1)$ of Definition \ref{AIR} is satisfied.

 Since ann-faith.dim$T\geq n$, there is a Hom$_{A}(-,T)$-exact exact sequence $A\stackrel{g_{0}}\to T_{0}\stackrel{g_{1}}\to T_{1}\to\cdots\stackrel{g_{n-1}}\to T_{n-1}$ where $T_{i}\in\mr{add}(T)$ for $0 \leq i \leq n-1$. Let $C=\mr{Coker}g_{n-1}$ and consider the exact sequence  $A\stackrel{g_{0}}\to T_{0}\stackrel{g_{1}}\to T_{1}\to\cdots\stackrel{g_{n-1}}\to T_{n-1}\to C\to 0$. By the definition of $\mr{Pres}^{n}(T)$, we see that $C\in\mr{Pres}^{n}(T)$ and hence $C\in\mr{Appres}^{\infty}(T)$ since $T$ is $n$-quasi-tilting. It follows that there is an exact sequence $0\to C^{'}\to T_{C}\to C\to 0$ with $C^{'}\in\mr{Appres}^{\infty}(T)$ and $T_{C}\in\mr{add}(T)$. We claim that the above exact sequence splits. In fact, by Lemma \ref{DA}$(2)$, for any $G\in\mr{Appres}^{n}(T)$, the induced sequence $\mr{Hom}_{A}(T_{n-1},G)\to\cdots\to \mr{Hom}_{A}(T_{0},G)\to\mr{Hom}_{A}(A,G)\to 0$ is exact. Obviously, this implies that Ext$^{1}_{A}(C,G)=0$. In particular, Ext$^{1}_{A}(C,C')=0$ and so the exact sequence $0\to C^{'}\to T_{C}\to C\to 0$ splits. It follows that $C\in\mr{add}(T)$. Now we see that the condition $(2)$ of Definition \ref{AIR} is satisfied. Thus $T$ is $n$-AIR-tilting.
\ \hfill $\Box$

\vskip 10pt

 $n$-Silting modules over any rings are introduced in \cite{MAO}. Let $T\in\mr{mod}A$ and  $\sigma_{_{T}}$ be an $(n+1)$-projective presentation of  $T$. Recall that $T$ is $ partial ~n\text{-}silting$ with respect to $\sigma_{_T}$ if $T\in\mc{D}_{\sigma_{_T}}$. Further, $T$ is $n\text{-}silting$ with respect to $\sigma_{_T}$ if $\mc{D}_{\sigma_{_{T}}}=\mr{Pres}^{n}(T)$. 
 
 By \cite[Lemma 3.5]{MAO}, if $T$ is $n$-silting, then $\mr{Pres}^{n}(T)=\mr{Pres}^{n+1}(T)$. This fact together with $\mc{D}_{\sigma_{_T}}\subseteq \mr{KerExt}_A^{1\le i\le n}(T,-)$ implies that $T$ is $n$-quasi-tilting.  In particular,  $\mc{D}_{\sigma_{_{T}}}=\mr{Appres}^{\infty}(T)$ if  $T$ is $n$-silting.

The following result provides a characterization of $n$-silting modules.

\bg{Pro}\label{n-silting1}
Let $T\in\mr{mod}A$ { and $\sigma_{_{T}}$ be an $(n+1)$-projective presentation of  $T$}. Then $T$ is an $n$-silting module  {with respect to $\sigma_{_T}$} if and only if $T\in\mc{D}_{\sigma_{_{T}}}$ and $\mc{D}_{\sigma_{_{T}}}\subseteq\mr{Gen}(T)$.
\ed{Pro}
\Pf. It is enough to show the sufficiency. As $\mc{D}_{\sigma_{_{T}}}\subseteq\mr{Gen}(T)$, for any $M\in \mc{D}_{\sigma_{_{T}}}$, there is a $\mr{Hom}(T,-)$-exact exact sequence $0\to M'\to T_M\to M\to 0$ with $T_M\in\mr{add}(T)$. It follows $M'\in \mc{D}_{\sigma_{_{T}}}$ too, by Lemma \ref{n-image} (5). Then one easily see that $\mc{D}_{\sigma_{_{T}}}\subseteq\mr{Pres}^n(T)$. As $\mr{Pres}^n(T)\subseteq\mc{D}_{\sigma_{_{T}}}$ already holds by Lemma \ref{n-image} (2), we obtain that $\mr{Pres}^n(T)=\mc{D}_{\sigma_{_{T}}}$ and hence $T$ is $n$-silting.
\ \hfill $\Box$

\vskip 10pt 
 We want to show the relations between $n$-AIR-tilting modules and $n$-silting  modules over an artin algebra.


\bg{Th}\label{main2}
Let $T\in\mr{mod}A$. Then the following are equivalent.

 $(1)$ $T$ is a strongly $n$-AIR-tilting module.

 $(2)$  $T$ is an  $n$-AIR-tilting module and satisfies $\mc{D}_{\sigma_{_{T}}}\subseteq \mc{D}_{\sigma_{_{A}}}$.
 
 $(3)$  $T$ is an $n$-silting module  and satisfies ann-faith.dim$T\geq n$.
\ed{Th}
\Pf. $(1)\Rightarrow (2)$. Clearly.
%

$(2)\Rightarrow (3)$  Since $T$ is an $n$-AIR-tilting module,  it holds $T\in\mc{D}_{\sigma_{_T}}$. By Lemma \ref{n-image}, we obtain Pres$^{n}(T)\subseteq\mc{D}_{\sigma_{_T}}$.
Now take any $M\in{\mc{D}}_{\sigma_{_T}}$, then there is an integer $m$ such that there exists an epimorphism $f
:A^{m}\to M$. Consider the following diagram, where the first row is the exact sequence given by $\sigma_A^m$.
$$
 \xymatrix{
&A^{m}\ar[r]^{g^{m}}\ar[r]\ar[d]^{f}& T^{m}_{0}\ar[r]\ar@{.>}[dl]^{h}     &  T_{1}^{m} \ar[r]            &  \cdots\ar[r] &T_{n}^{m} \ar[r]&0&&                      \\
 &M                     &  && & & && \\
}
$$
 By assumption,  $M\in{\mc{D}}_{\sigma_{_A}}$.  Hence there is a map $h:T^{m}_{0}\to M$ such that $f=hg^{m}$. It follows that $h$ is surjective. Then $\mc{D}_{\sigma_{_T}}\subseteq \mr{Gen}(T)$. By Lemma \ref{n-image}, ${\mc{D}}_{\sigma_{_T}}$ is closed under kernels of Hom$_{A}(T,-)$-epic epimorphisms. Then there is a Hom$_{A}(T,-)$-exact exact sequence $0\to M_{1}\to T_{0}^{'}\to M\to0 $ with $T_{0}^{'}\in\mr{add}(T)$ and $M_{1}\in\mc{D}_{\sigma_{_T}}$.
 Repeating this process to $M_{1}$ and so on, we obtain $\mc{D}_{\sigma_{_{_T}}}\subseteq \mr{Pres}^{n}(T)$. Thus, $T$ is $n$-silting.

$(3)\Rightarrow (1)$ Since $T$ is an $n$-silting module (with resp. to $\mc{D}_{\sigma_{_T}}$), we see that $T\in\mc{D}_{\sigma_{_T}}$. Since ann-faith.dim$T\geq n$, we obtain that there is a Hom$_{A}(-,T)$-exact exact sequence $A\to T_{0}\to T_{1}\to\cdots\to T_{n-1}$. Now the same method as we just use in the last part of the proof of Theorem \ref{AIR-snqt} shows that the above exact sequence gives the $\sigma_{_{A}}$ such that $T\in \mc{D}_{\sigma_{_A}}$.

Moreover, $T$ is an $n$-silting module implies that $\mc{D}_{\sigma_{_T}}=\mr{Pres}^{n}(T)=\mr{Appres}^{\infty}(T)$. And by Lemma \ref{DA}$(2)$, we obtain that
 $\mc{D}_{\sigma_{_A}}=\mr{Appres}^{\infty}(T)$. So together we see that $\mc{D}_{\sigma_{_T}}=\mc{D}_{\sigma_{_A}}$. Hence, $T$ is a strongly $n$-AIR-tilting module.
\ \hfill $\Box$


\vskip 10pt

By the proof of the second part of the above theorem, we obtain the following result which partly generalizes $\mr{\cite[Proposition~3.11]{AMV}}$.

\bg{Cor}\label{silting-AIR1}
Let $T\in\mr{mod}A$ { and $\sigma_{_{T}}$ be an $(n+1)$-projective presentation of  $T$}. Assume that $T$ be a partial  $n$-silting module {with respect to $\sigma_{_T}$}. Assume that there is a $\mr{Hom}_{A}(-,\mc{D}_{\sigma_{_T}})$-exact exact sequence $A\to T_{0}\to T_{1}\to\cdots\to T_{n}\to 0$, where $T_{i}\in\mr{add}(T)$ for $0\leq  i \leq n$. Then $T$ is an $n$-silting  module.
\ed{Cor}

{Recall that  $T\in\mr{mod}A$ is a $\ast^{n}$-{\it module} \cite{Wei4} if it satisfies the following conditions:}

{ $(1)$ $\mr{Pres}^{n}(T)=\mr{Pres}^{n+1}(T)$;}

{ $(2)$ For any exact sequence $0\to M\to T^{'}\to N\to 0$ in mod$A$, where  $M\in\mr{Pres}^{n}(T)$ and  $T^{'}\in\mr{add}(T)$, the induced sequence $0\to\mr{Hom}_{A}(T,M)\to \mr{Hom}_{A}(T,T^{'}) \to \mr{Hom}_{A}(T,N)\to 0$ is exact.}

It is clear that if $T$ is a $\ast^{n}$-module and let $\bar{A}=A/\mr{ann}(T)$, then $T_{\bar{A}}$ is also a $\ast^{n}$-module.

Let $T\in\mr{mod}A$ and $T$ be an $n$-quasi-tilting module. It is easy to see that $T$ is a $\ast^{n}$-module by the involved definitions. In particular,  $T_{\bar{A}}$ is a $\ast^{n}$-module, where $\bar{A}=A/\mr{ann}(T)$.



\bg{Pro}\label{nqt-tilting}
Let $T\in\mr{mod}A$ and $T$ be an $n$-quasi-tilting module.  Assume that there is a $\mr{Hom}_{A}(-,T)$-exact exact sequence $A\to T_{0}\to T_{1}\to \cdots\to T_{n}\to 0$, where $T_{i}\in\mr{add}(T)$ for $0 \leq i \leq n$. Then $T_{\bar{A}}$ is an $n$-tilting module, where $\bar{A}=A/\mr{ann}(T)$.
\ed{Pro}
\Pf. As $T_{\bar{A}}$ is already a $\ast^{n}$-module, we only need to show that injective $\bar{A}$-modules are contained  in Pres$^{n}_{\bar{A}}(T)$, by \cite{Wei4}. By  assumption, we have the Hom$_{\bar{A}}(-,T)$-exact sequence $\sigma^{'}_{\bar{A}}:0\to \bar{A}\to T_{0}\to T_{1}\to \cdots\to T_{n}\to 0$ in $\bar{A}$-modules. Let ${\mc{D}}_{\sigma^{'}_{\bar{A}}}:=\{M\in\mr{mod}{\bar{A}}| ~0\to \mr{Hom}_{\bar{A}}(T_{n},M)\to \mr{Hom}_{\bar{A}}(T_{n-1},M)\to\cdots\to\mr{Hom}_{\bar{A}}(\bar{A},M)\to 0~\mr{is~exact}\}$.  For any injective $\bar{A}$-module $I$, there is an integer $m$ such that $\bar{A}^{m}\to I$ is surjective.
Then we have the following commutative diagram
$$
 \xymatrix{
0\ar[r]&\bar{A}^{m}\ar[r]\ar[d]^{f}& T_{1}^{m} \ar@{.>}[dl]^{g}   &             &                      \\
 &I                     &  && \\
}
$$
Hence $g$ is surjective.  So,  $I\in\mr{Gen}T_{\bar{A}}$. Then there is a Hom$_{\bar{A}}(T,-)$-exact exact sequence $0\to K_{1}\stackrel{\alpha}\to T^{1}\stackrel{\beta}\to I\to 0$, where $T^{1}\in\mr{add}T_{\bar{A}}$.
  By Lemma \ref{DA}, we obtain $K_{1}\in\mc{D}_{\sigma^{'}_{\bar{A}}}$. Repeating the process to $K_{1}$ and so on, we obtain  $I\in\mr{Appres}^{n}(T_{\bar{A}})\subseteq \mr{Pres}^{n}(T_{\bar{A}})$. So,  $T$ is an $n$-tilting $\bar{A}$-module.
\ \hfill $\Box$

\vskip 10pt
We remark that, in case $n=1$, all of the notions, i.e., (strongly) 1-AIR-modules, (strongly) 1-quasi-tilting modules and 1-silting modules coincide with each other, for the artin algebra $A$ \cite{AMV,Wei3}. In particular, they are the same as support $\tau$-tilting modules.

\vskip 30pt
\section{Generalized two-term silting complexes}

In this section, we study the relation between $n$-AIR-tilting modules and silting complexes. As a corollary, we show that if the rank question for silting complexes has a positive answer, then strongly $n$-AIR-tilting modules and $n$-AIR-tilting modules coincide with each other, or equivalently, if strongly $n$-AIR-tilting modules and $n$-AIR-tilting modules are different, then the answer to the rank question for silting complexes is negative.

Recall that a complex $P^{\bullet}\in{D}^{b}(A)$ is $silting$ \cite{AI,KV,Wei2}if it satisfies the following conditions:

\begin{verse}
$(1)$ $P^{\bullet}$ is isomorphism to an object in ${K}^{b}(\mc{P}_A)$;

$(2)$ $\mr{Hom}_{{D}(A)}(P^{\bullet},P^{\bullet}[i])=0$ for all $i>0$;

$(3)$ $P^{\bullet}$ generates ${K}^{b}(\mc{P}_A)$, i.e., $\langle \mr{add}_{D(A)}{P^{\bullet}}\rangle={K}^{b}(\mc{P}_A)$, where $\langle\mr{add}_{D(A)}{P^{\bullet}}\rangle$  denotes the smallest triangulated subcategory containing add$P^{\bullet}$ in the  derived category.
\end{verse}

A complex $P^{\bullet}$ is said to be {\it presilting} if it is satisfies the above first two conditions.

A complex in $D(A)$ is said to satisfy the rank condition, if the number of it's distinct indecomposable direct summands is equal to the rank of the Grothendieck group $K_{0}(A)$, or equivalently, the numbers of distinct indecomposable direct summands of the complex and the regular module $A$ are the same. It is shown in \cite{AI} that silting complexes always satisfy the rank condition. Conversely, one asks {\it if a presilting complex that satisfies the rank condition is always silting?} This is called the rank question for silting complexes. The particular case, that is, the rank question for $n$-tilting modules was considered in \cite{RS}. 

A way to study the rank question is to investigate if a presilting complex can be completed to a silting complex, or equivalently, if it is a direct summand of a silting complex. Recall that a complex $P^{\bullet}\in {D}^{b}(A)$ is said to be {\it two-term} if it is isomorphic to a complex of the form $0\to P\to Q\to 0$ with $P,Q\in\mathcal{P}_A$, up to shifts. It is well known that every two-term presilting complex can be completed to a two-term silting complex, see for instance \cite{Wei2}. However, in the general case, it is not true and the first counterexample is given, very recently,  in \cite{LZ}, where it is shown that there is a presilting complex of length 3 (indeed, a pretilting module of projective dimension 2) which cannot be completed to any silting complex.

We remark that the rank question for silting complexes, even the particular case for general tilting modules, remains open in general.

We begin with the following lemma which characterizes the class $\mc{D}_{\sigma_{_T}}$ in the level of derived category.

\bg{Lem}\label{5-1}
Let $T\in\mr{mod}A$ and $\sigma_{_{T}}$ be an   $(n+1)$-projective presentation of $T$.

$(1)$ $\mc{D}_{\sigma_{_{T}}}={\sigma_{_{T}}^{\bullet}}^{\perp_{>0}}\cap\mr{mod}A$ and  $\mc{D}_{\sigma_{_{T}},k}={\sigma_{_{T}}^{\bullet}}^{\perp_{>n-k}}\cap\mr{mod}A$, for each $1\le k\le n$.

 $(2)$ For any $1\le k\le n$, an $A$-module $M$ belongs to $\mc{D}_{\sigma_{_{T}},k}$  if and only if for some $($respectively, all$)$  $(k+1)$-projective  presentation$($s$)$  $\omega$ of $M$ the the condition $\mr{Hom}_{D(A)}(\sigma_{_{T}}^{\bullet}, \omega^{\bullet}[i])=0$ for $i>n-k$ is satisfied.

\ed{Lem}
\Pf. (1) Note that $\sigma^{\bullet}: 0\to P_n\to\cdots\to P_0\to 0$ is in $K^b(\mathcal{P}_A)$, so any homomorphism in $D(A)$ starting from $\sigma^{\bullet}$ can be represented by a homomorphism in $K(A)$. Thus, for any $M\in \mr{mod}A$ and any integer $i$, one can see that a homomorphism $f: \sigma^{\bullet}\to M[i]$ is 0 in $D(A)$ if and only if $f$ is homotopic to 0. But the later is clearly equivalent to that the induced sequence $\mr{Hom}_A(P_{i-1},M)\to \mr{Hom}_A(P_{i},M)\to \mr{Hom}_A(P_{i+1},M) $ is exact, where $P_i=0$ if $i<0$ or $i>n$. The conclusion then follows from the involved definitions.

 (2) Assume that $\omega$ is of the form: $Q_{k}\stackrel{g_{k}}\to Q_{k-1}\to\cdots\stackrel{g_1}\to Q_{0}\stackrel{g_0}\to M\to 0$. Then $\omega^{\bullet}$ is the complex $0\to Q_{k}\stackrel{g_{k}}\to Q_{k-1}\to\cdots\stackrel{g_1}\to Q_{0}\to 0$.

The only-if part.  For any $i>n-k$, let $h^{\bullet}\in \mr{Hom}_{D(A)}(\sigma_T^{\bullet},\omega^{\bullet}[i]) \simeq \mr{Hom}_{K(A)}(\sigma_T^{\bullet},\omega^{\bullet}[i])$, and consider the following diagram. 
$$
 \xymatrix{
 &0\ar[r]&P_{n}\ar[r]^{f_{n}}\ar@{.>}[dl]_{s_{n-i+1}}\ar[d]|{h_{n-i}}& P_{n-1}\ar[r]^{f_{n-1}}\ar@{.>}[dl]|{s_{n-i}}\ar[d]|{h_{n-i-1}}&\cdots\ar[r]&P_{i+1}\ar[r]^{f_{i+1}}\ar[d]|{h_{1}} \ar@{.>}[dl]|{s_{i}}                &  P_{i} \ar[d]|{h_{0}}\ar[r]^{f_i}\ar@{.>}[dl]|{s_{1}}&  P_{i-1}\ar[r]\ar@{.>}[dl]^{s_0}\ar@{.>}[dd]^{h}&\cdots                       \\
 \cdots\ar[r]&Q_{n-i+1}\ar[r]&Q_{n-i}\ar[r]&\cdots\ar[r]&Q_{2}\ar[r]^{g_{2}}&Q_{1}\ar[r]^{g_{1} }                       & Q_{0}\ar[dr]_{g_0}\ar[r] &0\ \ \ \  &  \\
 &&&&&                     &  &  M&  \\
}
$$
Since $g_0 h_0 f_{i+1}=g_{0}g_{1}h_{1}=0 $ and $M\in\mc{D}_{\sigma_{_{T}},k}$, the induced sequence $\mr{Hom}_A(P_{i-1},M)\to \mr{Hom}_A(P_{i},M)\to \mr{Hom}_A(P_{i+1},M) $ is exact, and hence there is a map $h:P_{i-1}\to M$ such that $hf_{i}=g_{0}h_{0}$. Since $P_{i-1}$ is projective, there exists $s_{0}:P_{i-1}\to Q_{0}$ such that $h=g_{0}s_{0}$. Since $g_0(h_{0}-s_{0}f_{i})=0$, there is $s_{1}:P_{i}\to Q_{1}$ such that $h_{0}=g_{1}s_{1}+s_{0}f_{i}$. Apply the same process on $s_0$ to $s_1$, and so on,  we obtain $s_{j}: P_{j+i-1}\to Q_{j}$ such that $h_{j}=g_{j+1}s_{j+1}+s_{j}f_{j+i}$, where $2\leq j\leq n-i+1$. It follows that $h^{\bullet}=0$.

The if part. 
For any $n-k< i \leq n$, take $\alpha\in\mr{Hom}_{A}(P_{i},M)$ satisfying $\alpha f_{i+1}=0$, and consider the following diagram. $$
 \xymatrix{
 &0\ar[r]&P_{n}\ar[r]^{f_{n}}\ar@{.>}[dl]_{\alpha_{n-i+1}}&\cdots\ar[r]&P_{i+1}\ar[r]^{f_{i+1}} \ar@{.>}[dl]|{\alpha_{1}}                &  P_{i} \ar[d]|{\alpha}\ar[r]^{f_{i}}\ar@{.>}[dl]|{\alpha_{0}}\ar@{.>}[dll]^{s}&  P_{i-1}\ar[r]\ar@{.>}[dll]^{\ \ \ t}&\cdots                       \\
 \cdots\ar[r]&Q_{n-i+1}\ar[r]&\cdots\ar[r]&Q_{1}\ar[r]_{g_{1}}&Q_{0}\ar[r]_{g_{0} }                       & M &  &  \\
}
$$
Since $P_{i}$ is projective, there is a map $\alpha_{0}:P_{i}\to Q_{0}$ such that $\alpha=g_{0}\alpha_{0}$. Inductively, we obtain a chain map $\alpha^{\bullet}:\sigma_T^{\bullet}\to \omega^{\bullet}[i]$. By assumption, $\alpha^{\bullet}=0$, and hence there are $s:P_{i}\to Q_{1}$ and $t:P_{i-1}\to Q_{0}$ such that $\alpha_{0}=g_{1}s+tf_{i}$. Hence $\alpha=g_{0}\alpha_{0}=g_{0}tf_{i}$. This show that the induced sequence  $\mr{Hom}_A(P_{i-1},M)\to \mr{Hom}_A(P_{i},M)\to \mr{Hom}_A(P_{i+1},M) $ is exact. It follows that  $M\in\mc{D}_{\sigma_{_{T}},k}$. 
\ \hfill $\Box$

\vskip 10pt
We say that a module $M$ is {\it $n$-pre-AIR-tilting} if $T\in\mc{D}_{\sigma_{_{T}}}$. The following immediate corollary provides a characterization of $n$-pre-AIR-tilting modules.

\bg{Cor}\label{preAIR}
Let $M\in \mr{mod}A$. Then $M$ is $n$-pre-AIR-tilting if and only if $\omega^{\bullet}$ is presilting, where $\omega^{\bullet}$ be the associated truncated complex of the minimal $($resp. some$)$ $(n+1)$-projective presentation of $M$.
\ed{Cor}

In \cite{AIR}, the authors shows that
support $\tau$-tilting modules are in bijection with two-term silting complexes. Since support $\tau$-tilting modules are just $1$-AIR-tilting modules, it is natural to consider the relation between $n$-AIR-tilting modules and general silting complexes.

As a generalization of two-term complexes, we say that a complex $P^{\bullet}\in {D}^{b}(A)$ is  {\it generalized two-term} if it is isomorphic, up to shifts, to a complex of the form $0\to P_{n}\to P_{n-1}\to \cdots\to P_{0}\to 0$, where each $P_i\in\mathcal{P}_A$, such that  $\mr{H}^{-i}(P^{\bullet})=0$, for $0<i<n$.


The following result shows that every generalized two-term silting complex gives a strongly $n$-AIR-tilting module.

\bg{Pro}\label{silt-AIR}
Let $P^{\bullet}$ be a generalized two-term silting complex of the form $0\to P_n\to \cdots\to P_0\to 0$, up to shifts. Then $T=\mr{H}^{0}(P^{\bullet})$ is a strongly $n$-AIR-tilting module.
 \ed{Pro}
 
\Pf. Denote $\sigma_{_{T}}$ the obvious exact sequence $P_n\to \cdots\to P_0\to T\to 0$. Then $T$ is $n$-pre-AIR-tilting by Corollary \ref{preAIR}, since $\sigma_T^{\bullet}=P^{\bullet}$ is presilting.

We now show that there is an exact sequence $\sigma_A: A\to T_0\to \cdots\to T_n\to 0$ with $T\in \mc{D}_{\sigma_{_A}}$.

Since $P^{\bullet}$ is generalized $2$-term silting,
 there are   triangles $X^{\bullet}_{i}\stackrel{h_{i}^{\bullet}}\to P_{i}^{\bullet}\stackrel{g^{\bullet}_{i+1}}\to X_{i+1}^{\bullet}\to$ with each $P^{\bullet}_{i}\in\mr{add}_{D(A)} P^{\bullet}$ for $0\leq i \leq n$, where $X^{\bullet}_{0}=A$ and $X^{\bullet}_{n+1}=0$ by $\mr{\cite[Theorem~3.5]{Wei2}}$. Note that $X^{\bullet}_{n}=P^{\bullet}_{n}\in \mr{add}_{D(A)}P^{\bullet}$.
Let $f_0^{\bullet}=h^{\bullet}_{0}$ and $f_{i}^{\bullet}=h^{\bullet}_{i}g^{\bullet}_{i}$ for $1 \leq i \leq n$.

Applying the functor $\mr{Hom}_{K(A)}(A,-)$ to these triangles,
we claim that the induced sequence $$\mr{Hom}_{K(A)}(A,A)\stackrel{(A,f^{\bullet}_{0})}\longrightarrow \mr{Hom}_{K(A)}(A,P^{\bullet}_{0})\stackrel{(A,f^{\bullet}_{1})}\longrightarrow \mr{Hom}_{K(A)}(A,P^{\bullet}_{1})\to\cdots\to\mr{Hom}_{K(A)}(A,P^{\bullet}_{n})\to 0 $$
is exact.

Consider the following diagram
$$
 \xymatrix{
 && A\ar[d]^{f^{\bullet}} \ar@{.>}[dl]_{h} \ar@/_/@{.>}[ddd]_{g^{\bullet}  }          &  &  & & && &&  &\\
 &A\ar[r]^{f_{0}^{\bullet}} &P_{0}^{\bullet}\ar[rr]^{f^{\bullet}_{1}}\ar[dr]_{g^{\bullet}_{1}}          &    &  P_{1}^{\bullet} \ar[dr]_{g^{\bullet}_{2}}\ar[rr]^{f^{\bullet}_{2}}&&  P^{\bullet}_{2 }\ar[r] &\cdots\ar[r]& P_{n-1}^{\bullet}\ar[rr]^{f_{n}}\ar[dr]^{g^{\bullet}_{n}}&&P^{\bullet}_{n} &                      \\
 &&                  &X^{\bullet}_{1}\ar[ur]_{h_{1}^{\bullet} }  &  & X^{\bullet}_{2}\ar[ur]_{h^{\bullet}_{2}}& &X^{\bullet}_{n-1}\ar[ur]_{h_{n-1}}& &X^{\bullet}_{n}\ar[ur]^{\simeq}&  & \\
 &&  X^{\bullet}_{2}[-1]\ar[ur]_{\alpha^{\bullet}}                &  &  & & && &&  &}
$$

Firstly, for any $f^{\bullet}\in\mr{Hom}_{K(A)}(A,P^{\bullet}_{0})$ such that $f_{1}^{\bullet}f^{\bullet}=0$, that is $h^{\bullet}_{1}g^{\bullet}_{1}f^{\bullet}=0$, there exists $g^{\bullet}:A\to X_{2}^{\bullet}[-1]$ such that $g^{\bullet}_{1}f^{\bullet}=\alpha^{\bullet} g^{\bullet}$. Applying Hom$_{K(A)}(A,- )$ to the triangle $X^{\bullet}_{2}\to P^{\bullet}_{2}\to X^{\bullet}_{3}\to$, we obtain the induced exact sequence $$\mr{Hom}_{K(A)}(A,P^{\bullet}_{2}[-2])\to\mr{Hom}_{K(A)}(A,X^{\bullet}_{3}[-2])\to \mr{Hom}_{K(A)}(A,X_{2}^{\bullet}[-1])\to  \mr{Hom}_{K(A)}(A,P^{\bullet}_{2}[-1]).$$
Since $\mr{Hom}_{K(A)}(A,P^{\bullet}_{2}[-2])=\mr{Hom}_{K(A)}(A,P^{\bullet}_{2}[-1])=0$ by the assumption, we obtain that $\mr{Hom}_{K(A)}(A,X^{\bullet}_{3}[-2])\simeq \mr{Hom}_{K(A)}(A,X_{2}^{\bullet}[-1])$. Inductively, by applying $\mr{Hom}_{K(A)}(A,- )$ to the triangles $X^{\bullet}_{i}\to P^{\bullet}_{i}\to X^{\bullet}_{i+1}\to$, we obtain that
$$\mr{Hom}_{K(A)}(A,X_{2}^{\bullet}[-1])\simeq \mr{Hom}_{K(A)}(A,X_{3}^{\bullet}[-2]) \simeq \cdots\simeq\mr{Hom}_{K(A)}(A,X^{\bullet}_{n-1}[-(n-2)])=0,$$
where the last equality holds because there is an induced exact sequence
{\footnotesize{$$0=\mr{Hom}_{K(A)}(A,P^{\bullet}_{n}[-(n-1)])\to \mr{Hom}_{K(A)}(A,X^{\bullet}_{n-1}[-(n-2)])\to\mr{Hom}_{K(A)}(A,P^{\bullet}_{n-1}[-(n-2)])=0.$$}}
\noindent Then $g^{\bullet}=0=g^{\bullet}_{1}f^{\bullet}$, and  hence there is some $h:A\to A$ such that $f^{\bullet}=f_{0}^{\bullet}h$. It follows that  $\mr{Hom}_{K(A)}(A,A)\to \mr{Hom}_{K(A)}(A,P^{\bullet}_{0})\to \mr{Hom}_{A}(A,P^{\bullet}_{1})$ is exact.
Similarly,  we could prove that Hom$_{K(A)}(A,P^{\bullet}_{i-1})\to\mr{Hom}_{K(A)}(A,P^{\bullet}_{i})\to\mr{Hom}_{K(A)}(A,P^{\bullet}_{i+1})$ is exact, since  H$^{-j}(P^{\bullet})=0$  for $1\leq i \leq n-1$ and $1\leq j \leq n-1$.

Since $\mr{H}^{i}(A)=\mr{H}^{i}(P^{\bullet}_{j})=0$, we obtain $\mr{H}^{i}(X^{\bullet}_{j})=0$ for all $i> 0$ and $0\leq j \leq n $. Then $\mr{Hom}_{K(A)}(A,X^{\bullet}_{n-1}[1])=0$; hence  $\mr{Hom}_{K(A)}(A, P^{\bullet}_{n-1})\to \mr{Hom}_{K(A)}(A,P^{\bullet}_{n})\to 0$ is exact.

 So, we obtain that the sequence  $\mr{Hom}_{K(A)}(A,A)\to \mr{Hom}_{K(A)}(A,P^{\bullet}_{0})\to \mr{Hom}_{K(A)}(A,P^{\bullet}_{1})\to\mr{Hom}_{K(A)}(A,P^{\bullet}_{n})\to 0 $   is exact. Since $\mr{Hom}_{K(A)}(A,P^{\bullet}_{i})\simeq\mr{H}^{0}(P^{\bullet}_{i})\in\mr{add}(T)$, we have an exact sequence $A\stackrel{g_{0}}\to T_{0}\stackrel{g_1}\to T_{1}\to\cdots\to T_{n}\to 0~~(\ast)$, where $T_{i}\in\mr{add}(T)$ for $0\leq i\leq n$. 
 
 Now, we prove that the exact sequence  $(\ast)$ is $\mr{Hom}_{A}(-,T)$-exact. Since $T$ is $n$-pre-AIR-tilting, $T\in\mr{KerExt}_A^{1\le i\le n}(T,-)$. Then, it is enough to show that the induced sequence $\mr{Hom}_{A}(T_{0},T)\to \mr{Hom}_{A}(A,T)\to 0$ is exact.

Take any $f\in\mr{Hom}_{A}(A,T)$. Then $f$ induces a chain map $g^{\bullet}: A\to P^{\bullet}$, as follows.
$$
 \xymatrix{
 &&     &  A \ar[d]_{f}\ar[r]^{f^{\bullet}_{0}}  \ar@{.>}[dl]_{g}              &  P^{\bullet}_{0}\ar[r] &P^{\bullet}_{1}\ar[r] &\cdots\ar[r]&P^{\bullet}_{n}&                      \\
 \cdots\ar[r]&P_1\ar[r]                     & P_{0}\ar[r]_{f_{0}} &  T\ar[r]&0& && \\
}
$$
Note that  $\mr{Hom}_{K(A)}(X^{\bullet}_{k},P^{\bullet}[j])=0$ for $1\leq k \leq n$ and $j>0$, 
since  $\mr{Hom}_{K(A)}(A,P^{\bullet}[j])=0$ and $\mr{Hom}_{K(A)}(P^{\bullet}_{i},P^{\bullet}[j])=0$ for each $i$. In particular, $\mr{Hom}_{K(A)}(X^{\bullet}_{1}[-1],P^{\bullet})=0$. Then there exists $h^{\bullet}:P^{\bullet}_{0}\to P^{\bullet}$ such that $g^{\bullet}=h^{\bullet}f^{\bullet}_{0}$.

$$
 \xymatrix{
 & X_{1}^{\bullet}[-1]\ar[r] &A\ar[r]^{f^{\bullet}_{0}} \ar[d]^{g^{\bullet}}               &  P^{\bullet}_{0}\ar@{.>}[dl]^{h^{\bullet}}  \ar[r]                & &                       \\
 &                       &P^{\bullet}  &  & \\
}
$$
Applying H$^{0}(-)$ to the above diagram, we obtain the following commutative diagram

$$
 \xymatrix{
 & &A\ar[r]^{g_{0}} \ar[d]^{f}               &  T_{0}\ar@{.>}[dl]^{h}                  & &                       \\
 &                       &T  &  & \\
}
$$
So, we obtain that induced sequence $\mr{Hom}_{A}(T_{0},T)\to \mr{Hom}_{A}(A,T)\to 0$ is exact.
It follows that $(\ast)$ is the desired sequence $\sigma_{_A}$. In particular, we have that {ann-faith.dim}$T\ge n$.

By Theorem \ref{main2} and Proposition \ref{n-silting1} and the involved definitions, it remains to show that $\mc{D}_{\sigma_{_T}}\subseteq \mr{Gen}(T)$.

Since $\sigma_T^{\bullet}=P^{\bullet}: 0\to P_n\to \cdots\to P_0\to 0$ is a silting complex, we obtain that $P^{\bullet \bot_{ >0}} = \mr{Pres}^{n}_{D^{\leq0}}(P^{\bullet}) :=\{X^{\bullet}\in\mr{D}(A)~|~\mr{there ~are ~triangles} ~X^{\bullet}_{i+1}\to P^{\bullet}_{i}\to X^{\bullet}_{i} \to, \mr{where}~P^{\bullet}_{i}\in\mr{add}P^{\bullet},X^{\bullet}_{n+1}\in{D}^{\leq0} ~\mr{and}~X^{\bullet}_{0}=X^{\bullet} ~\mr{for}~0\leq  i \leq n+1 \}$ by $\mr{\cite[Theorem~4.4]{Wei2}}$, where $D^{\le 0}$ is the subcategory of complexes with homologies vanishing for all positive degrees.
Take any $X\in \mc{D}_{\sigma_{_T}}$. Then $X\in P^{\bullet \bot_{>0} }\cap\mr{mod}A$ by Lemma \ref{5-1}. Thus, the above arguments implies that there is a triangle $X^{\bullet}_{1}\to P^{\bullet}_{1}\to X\to$ with $X^{\bullet}_{1}\in{D}^{\leq0}$ and $P^{\bullet}_{1}\in\mr{add}P^{\bullet}$. Applying H$^{0}(-)$ to this triangle, we obtain an exact sequence H$^{0}({P_{1}^{\bullet}})\to X\to 0$. Since  H$^{0}({P^{\bullet}_{1}})\in\mr{add}(T)$, we see that $X\in\mr{Gen}(T)$. It follows that $\mc{D}_{\sigma_{_T}}\subseteq \mr{Gen}(T)$.
\ \hfill $\Box$

\vskip 10pt
Conversely, one want to know if a strongly $n$-AIR-tilting module lifts a silting complex. We have the following partial answer to this question.

\bg{Pro}\label{AIR-2-silt}
Let $T\in\mr{mod}A$ and $P^{\bullet}$ be the associated truncated complex of the minimal $(n+1)$-projective presentation of $T$. Let $T$ be an $n$-AIR-tilting module. Assume that the rank question for silting complexes has a positive answer, then there is a projective module $Q$ such that $P^{\bullet} \oplus Q[n]$ is a generalized $2$-term silting complex.
\ed{Pro}
\Pf. Denote $|M|$ to be the number of distinct indecomposable direct summands of the object $M$. By Proposition \ref{AIR-tilting}, $T$ is a tilting $\bar{A}$-module, where $\bar{A}=A/\mr{ann}(T)$. Then $|\bar{A}|=|T|$. Let $Q$ be the maximal basic projective module such that $\mr{Hom}_A(Q,T)=0$ and let $Q=eA$ for some idempotent $e$. Then $|\bar{A}|=|A/\langle e\rangle|=|A|-|Q|$. It follows that $|A| = |Q|+|T|$. By the assumption we know that $|T|=|P^{\bullet}|$. Hence we have that $|A| = |Q[n]|+|P^{\bullet}| = |Q[n]\oplus P^{\bullet}|$.

As $T$ is an $n$-AIR-tilting module,  $P^{\bullet}$ is a presilting complex by Corollary \ref{preAIR}. By the construction, $Q$ is a projective module such that $\mr{Hom}_A(Q,T)=0$, so one can easily check that $P^{\bullet} \oplus Q[n]$ is also presilting. Then we conclude that $P^{\bullet} \oplus Q[n]$ is a presilting complex satisfying the rank condition. It follows that  $P^{\bullet}\oplus Q[n]$ is a silting complex (which is clearly generalized two-term) by our assumptions.
\ \hfill $\Box$

\vskip 10pt
Two objects $M$ and $N$ in an additive category are said to be equivalent if $\mr{add}(M)=\mr{add}(N)$.
We combining the above two results to the following theorem, which is a generalization of \cite[Theorem 3.2]{AIR}.

\bg{Th} \label{2-silt-AIR-th} Let $n$ be a fixed positive integer. Let {\bf g2-silt$_n$(A)} be the class of generalized two-term silting complexes of the form $0\to P_n\to\cdots\to P_0\to 0$ with each $P_i\in \mc{P}_A$, up to shifts and equivalences, and let {\bf sAIR-tilt$_n$(A)} $(${\bf AIR-tilt$_n$(A)}, resp.$)$ be the class of strongly $n$-AIR-tilting modules $($$n$-AIR-tilting modules, resp.$)$, up to equivalences.  Then there are the following maps, where $\Phi$ is given in Proposition \ref{silt-AIR} and $\Psi$ is given in Proposition \ref{AIR-2-silt} with the additional assumption on the latter map that the rank question for silting complexes has a positive answer. In the case, the morphisms $\Phi$ and $\Psi$ are inverse maps.

\vskip 8pt
\centerline{{\bf g2-silt}$_n(A)\stackrel{\Phi}{\to}$ {\bf sAIR-tilt}$_n(A)\subseteq$ {\bf AIR-tilt}$_n(A)\stackrel{\Psi}{\to}$ {\bf g2-silt}$_n(A)$}
\ed{Th}

\Pf. It remains to show that $\Phi$ and $\Psi$ are inverse maps under the assumptions, but this follows from the constructions of the two maps.
\hfill$\Box$


\vskip 20pt
\section{Examples and questions}

The first example coming from \cite{Wei5}.

\bg{Exm} \rm{Every simple module which is pretilting is strongly $n$-AIR-tilting.}
\ed{Exm}

\Pf. By the proof of \cite[Theorem 1]{Wei5}, such a simple module is the last homology of a generalized two-term silting complexes. Hence, it is strongly $n$-AIR-tilting by Proposition \ref{silt-AIR}.
\hfill$\Box$

\vskip 10pt
Indeed, by the same reason, modules presented in \cite[Theorem 2]{Wei5} are also strongly $n$-AIR-tilting.

The module in the following example is used in \cite{LZ} as a counterexample to the complement question for presilting complexes. We see that the module can not be $n$-AIR-tilting, but it is strongly $n$-quasi-tilting.

\bg{Exm} \rm{Let $A = kQ/I$, where
\[Q = \xymatrix{1 \ar@/^0.5pc/[r]^{x_1} \ar@/_0.5pc/[r]_{y_1} & 2 \ar@/^0.5pc/[r]^{x_2} \ar@/_0.5pc/[r]_{y_2} & 3 },\ \text{and}\
I = \langle x_1x_2, y_1y_2 \rangle. \]
Let $T$ be the indecomposable module represented by 
$\xymatrix{1\ar[r]^{x_1}&2\ar[r]^{y_2}&3}$. Then 

(1) $T$ is a strongly  $n$-quasi-tilting module for any $n\ge 2$.

(2) $T$ is $n$-pre-AIR-tilting for any $n\ge 2$.

(3) $T$ is not $n$-silting for any $n$.

(4) $T$ is not $n$-AIR-tilting  for any $n$.
}
\ed{Exm}

\Pf. The module $T$ has the minimal projective resolution as follows, where the image of the morphism $f$ is $\xymatrix{2\ar[r]^{x_2}&3}$.

\centerline{$\xymatrix{
0\ar[r]&~3~\ar[r]&~(2\ar@<.5ex>[r]^{x_2}\ar@<-.5ex>[r]_{y_2}&{^3_3})~\ar[r]^f&~(1\ar@<.5ex>[r]^{x_1}\ar@<-.5ex>[r]_{y_1}&{^2_2}\ar@<.5ex>[r]^{y_2}\ar@<-.5ex>[r]_{x_2}&{^3_3})~\ar[r]&~(1\ar[r]^{x_1}&2\ar[r]^{y_2}&3)~\ar[r]&~0}$}

As shown in \cite{LZ}, $T$ is pretilting with the projective dimension 2. Thus, (2) holds by the involved definition. 

Note that $\mr{Gen}(T) = \mr{add}(1\oplus (\xymatrix{1\ar[r]^{x_1}&2}) \oplus T)$ and that $\mr{Pres}(T) = \mr{add}(T) = \mr{Pres}^n(T)$ for any $n\ge 2$, so we get that $T$ is $n$-quasi-tilting for any $n\ge 2$. Since $\mr{Hom}_A(3,1\oplus (\xymatrix{1\ar[r]^{x_1}&2})) =0$, one easily see that $\mr{Gen}(T)\subseteq \mr{KerExt}_A^2(T,-)$. Thus $T$ is strongly $2$-quasi-tilting. It is also easy to see that $T$ is strongly $n$-quasi-tilting for any $n\ge 3$. So (1) follows.

Since, for any $n\ge 2$, every indecomposable injective module is  clearly contained in $\mc{D}_{\sigma_{_T}}$, where $\sigma_T$ is the associated truncated complex of the minimal projective resolution of length $n+1$, and since that $\mr{Gen}(T)$ doesn't contain all injective modules, we get that $T$ is not $n$-silting for any $n\ge 2$. In case that $n=1$, as the obvious  morphism $g: \xymatrix{(2\ar@<.5ex>[r]^{x_2}\ar@<-.5ex>[r]_{y_2}&{^3_3})~\ar[r]&~T}$ doesn't factor through the morphism $f$, $T$ is not 1-pre-AIR-tilting. Thus, $T$ is not 1-silting module too. Then the conclusion (3) follows.

Finally, let $h: A\to T_A$ be the minimal left $\mr{add}(T)$-approximation. Then one can check that $\mr{Coker}h=1\oplus (\xymatrix{1\ar[r]^{x_1}&2})$ and that $\mr{Hom}_A(\mr{Coker}h,T)=0$. So there is no exact sequence $A\to T_0\to \cdots\to T_n\to 0$ with each $T_i\in\mr{add}(T)$ for any $n\ge 1$. Thus (4) follows.
\hfill$\Box$

\vskip 10pt
The following example shows that $n$-quasi-tilting modules and strongly $n$-quasi-tilting modules are not the same.

\bg{Exm} \rm{Let $A = kQ/I$, where
\[Q = \xymatrix{1 \ar@<.5ex>[r]^{\alpha}  & 2 \ar@<.5ex>[r]^{\beta}  & 3\ar@/^0.5pc/[ll]^{\gamma}},\ \text{and}\
I = \langle \beta\gamma, \gamma\alpha \rangle. \]
Let $T$ be the indecomposable module represented by 
$\xymatrix{1\ar[r]^{\alpha}&2}$. Then $T$ is $2$-quasi-tilting but not strongly $2$-quasi-tilting.
}
\ed{Exm}

\Pf. The module $T$ has the minimal projective presentation of length 3 as follows,  where the image of the morphism $f$ is the simple module 1.

\centerline{$\xymatrix{
(1\ar[r]^{\alpha}&2\ar[r]^{\beta}&3)~\ar[r]^f&~(3\ar[r]^{\gamma}&1)~\ar[r]&~(1\ar[r]^{\alpha}&2\ar[r]^{\beta}&3)~\ar[r]&~(1\ar[r]^{\alpha}&2)~\ar[r]&~0}$}

It is easy to check that $\mr{Pres}(T) = \mr{add}(T) \subseteq \mr{KerExt}_A^{1\le i\le 2}(T,-)$, thus one can get that $T$ is $2$-quasi-tilting. But the simple module $1\in \mr{Gen}(T)$ is not in $\mr{KerExt}_A^2(T,-)$ since it is just the second syzygy of $T$, it follows that $T$ is not strongly $2$-quasi-tilting.
\hfill$\Box$

\vskip 10pt

At last, we list two questions remained in the paper.

\vskip 10pt

{\bf Question 1} {\it Are $n$-AIR-tilting modules and strongly $n$-AIR-tilting modules the same?}

\vskip 10pt
If Question 1 has a negative answer, then one can get a counterexample to the rank question for silting complexes, by Theorem \ref{2-silt-AIR-th}.

\vskip 10pt
{\bf Question 2} {\it Are $n$-silting modules also $n$-AIR-tilting modules?}

\vskip 10pt
If Question 2 has the positive answer, then Question 1 also have a positive answer, by Theorem \ref{main2}.

%
%
%
%

{\small

}

\end{document}